\documentclass{llncs}

\usepackage{hyperref}

\usepackage{cite}
\usepackage{graphicx}
\usepackage{framed}
\usepackage{amsfonts}
\usepackage{amsmath}
\usepackage{color}
\usepackage{algorithm}
\usepackage[noend]{algpseudocode}
\definecolor{mygreen}{rgb}{0,0.2,0}
\definecolor{mygray}{rgb}{0.95,0.95,0.95}
\definecolor{mymauve}{rgb}{0.58,0,0.82}

\definecolor{lbcolor}{rgb}{0.95,0.95,0.95}

\usepackage{listings,newtxtt}
\begin{document}

\title{Teaching fractals for gifted learners at age 12 \\by using novel technologies}%\\2019-09-04 --- Draft 7}
\author{Zolt\'an Kov\'acs\orcidID{0000-0003-2512-5793}}
\institute{
The Private University College of Education of the Diocese of Linz\\
Salesianumweg 3, A-4020 Linz, Austria\\
\email{zoltan@geogebra.org}
}

\maketitle              % typeset the title of the contribution

\begin{abstract}
A summary of an experimental course on fractals is given
that was held for young learners at age 12. The course was a part
of Epsilon camp, a program designed for very gifted students who have
already demonstrated high interest in studying mathematics. Prerequisites for
the course were mastery of Algebra I, experience and fluency in skills like
exponentials and square roots, solving equations.
Also, at least two preliminary years were required in a prior
Epsilon camp. The summary gives an overview of the flow of teaching,
the achieved results and some evaluation of the given feedback.

\keywords{fractals, XaoS, GeoGebra, 3D printing, basic calculus}

\end{abstract}
\section{Introduction}

Epsilon Camp (\url{https://epsiloncamp.org}) is an annual meeting of young learners in the United States that has been held
at different campus locations since 2011 for two weeks during summer vacation.
It was founded by George Reuben Thomas
who is currently the executive director of the camp. The camp program is organized
by an academic director who follows the mission statement of the camp. In a nutshell,
the young learners are divided into groups that reflect their age and preliminary skills.

Campers are placed into groups of 10--12 campers with
a dedicated counselor, who is usually a university student in mathematics. Camper groups within each level do not imply
any difference in ability. Camper levels are currently as follows:

\begin{enumerate}
        \item Pythagoras Level: 7- and 8-year-olds.
        \item Euclid Level: 9- to 11-year-olds who are new to camp or 9- to 11-year-olds who have previously completed any part of Pythagoras Level.
        \item Gauss Level: 10- to 12-year-olds who have completed the Euclid Level.
        \item Conway Level: 11- to 12-year-olds who have completed the Gauss Level.
\end{enumerate}

Typical curriculum topics for campers have included introductory topics to advanced mathematics,
including number theory, methods of proofs, voting theory, set theory, Euclidean
geometry, projective geometry, or hyperbolic geometry, among others. These topics
can vary from year to year, but the basic concept is that for higher levels some
non-introductory topics are also included.

In this paper we sketch up the flow of a course given at Epsilon Camp 2019,
introduced for the Conway Level, held 14--28 July 2019 at University of Colorado, Colorado Springs.

\section{Course description}
Eight students attended the course: 3 girls and 5 boys. The youngest learner was 11 and
the oldest 12 years old.

The course consisted of 10 classes, 80 minutes each. The students were taught
in a university classroom with whiteboard and projector access. The students
worked only with paper and pencil. Using calculators for the students in the camp was discouraged,
however in this course at some points it was still allowed.

The basic idea of the course was to explain some basic concepts of fractals,
with an emphasis on self-similarity. In mathematical sense some basic concepts from calculus
were used frequently. The dragon curve, the Koch curve,
the Koch snowflake, the Sierpi\'nski carpet, Cantor's ternary set and the Sierpi\'nski triangle
were introduced first. Then escape time fractals like the Mandelbrot set were discussed by defining
bifurcation via accumulation points.
After introducing complex value root computation, the cubic Newton fractal was shown and investigated
in detail.

Novel teaching methods like using computer algebra, a recent version of the \textit{XaoS} software tool \cite{xaos-gh}, or
3D printing, were also utilized in the teaching process.

\section{The used didactic method}
Epsilon Camp's mission statement discouraged teaching
in a definition--theorem--proof style. Instead, a challenge problem (namely, ``Can several copies of the dragon curve
tile the plane'') was introduced, and the underlying theory to study fractals
was built with extensive help of the students. The students' natural curiosity therefore led to
natural discovery of new knowledge. The learning process was facilitated through
the Socratic method by the teacher: many questions were raised in order to
follow the expected line of reasoning.

Homework helped to deepen knowledge or to prepare next day's activity. Since each
group of children was learning three different subjects and the available time
for doing homework was 70 minutes per day, it was not expected that students
would spend more than 20 minutes to solve all homework exercises. On the other hand,
some students used some extra time to finalize their assignments and turned them
in---sometimes with standards of high quality. %% see Josh's homeworks

Most assignments were purely mathematical; some were about to do some experiments
manually with paper folding or by using 3D printed objects.

\section{The course in detail}

In this section a detailed description of the course follows. A daily step-by-step
explanation of the topic is presented. 

\subsection*{Day 1}

As an introductory example a few iterations of the dragon curve were shown, by using colored paper stripes.
There was an immediate technical problem with the stripes that was identified early
enough (before the first day): the stripes being used were not robust enough to stand firmly
and therefore the right angles were inaccurate. However, as homework, the students
were asked to try to fold better instances of the first iterations. Unfortunately,
this activity was maybe not interesting enough for most of the students, only one participant
was working on this task---but, in fact, one of her interests during the camp was paper folding in general.

\begin{figure}
\begin{center}
\includegraphics[height=0.3\textwidth]{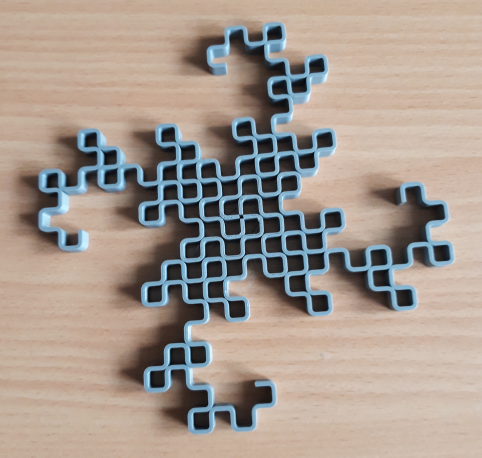}
\includegraphics[height=0.3\textwidth]{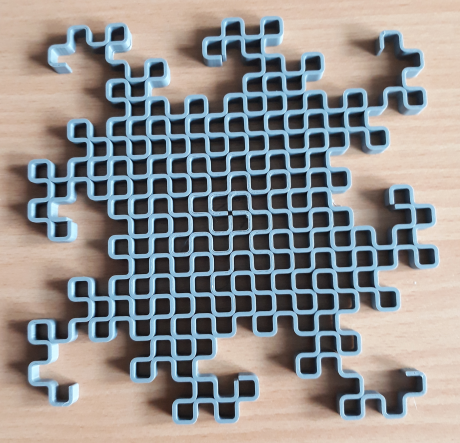}
\includegraphics[height=0.3\textwidth]{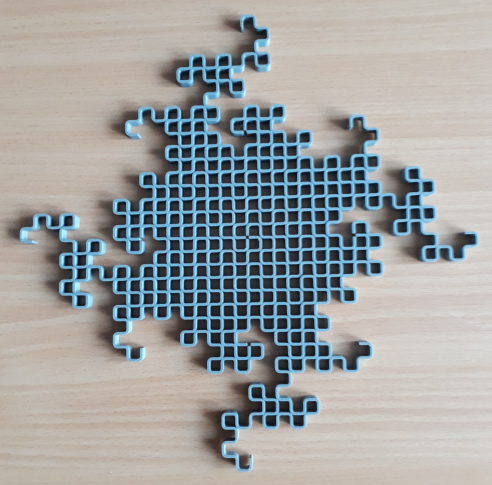}
\caption{An assembly of 4, 8 and 12 pieces of the 7th iteration of the dragon curve}
\label{3dragons}
\end{center}
\end{figure}

\begin{figure}
\begin{center}
\includegraphics[height=0.3\textwidth]{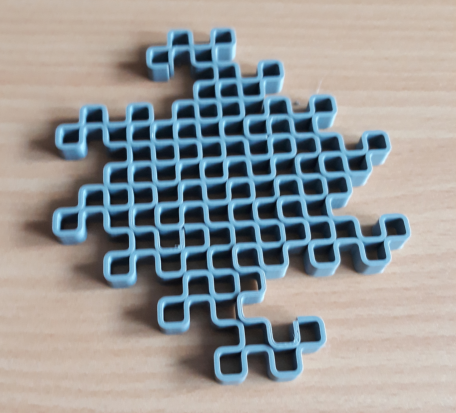}
\caption{An alternative assembly of 4 pieces of the 7th iteration of the dragon curve}
\label{dragon-d}
\end{center}
\end{figure}

One of the final aims of the course was to mention the plane filling property of the dragon curve.
To achieve this aim the author brought 15 pieces of 3D prints of the 7th iteration of the dragon curve and shared them
among the students. As a second homework for the group, students had to verify the plane filling property by
using their own pieces in a bigger picture. Fig.~\ref{3dragons} shows an example assembly of 4, 8 and 12 pieces
of prints, respectively. (Note that there is a high variety of possible assemblies that fill the plane.
For example, Fig.~\ref{dragon-d} shows a variant that contains 4 pieces.)

\begin{figure}
\begin{center}
\includegraphics[height=0.3\textwidth]{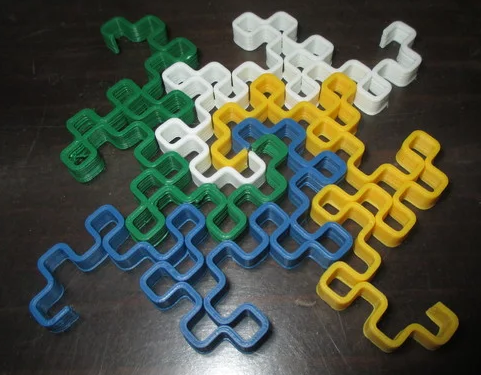}
\caption{Yet another alternative assembly of 4 pieces of the 7th iteration of the dragon curve---with different colors}
\label{thingiverse-color}
\end{center}
\end{figure}

The 7th iteration of the dragon curve can be downloaded
as a freely available STL file from \url{https://www.thingiverse.com/thing:1937697},
made by a contributor of the University of Maine at Farmington.
The 3D prints were proven to be interesting and robust enough. For a future version of the course
a mucher higher amount of prints could also be used. (One piece of print took 41 minutes on a HicTop CR-10
3D printer.) It would be useful to provide the prints in various colors (see Fig.~\ref{thingiverse-color}
for a photo provided by the contributor---it shows another variant of 4 pieces).

\begin{figure}
\begin{center}
\includegraphics[width=0.4\textwidth]{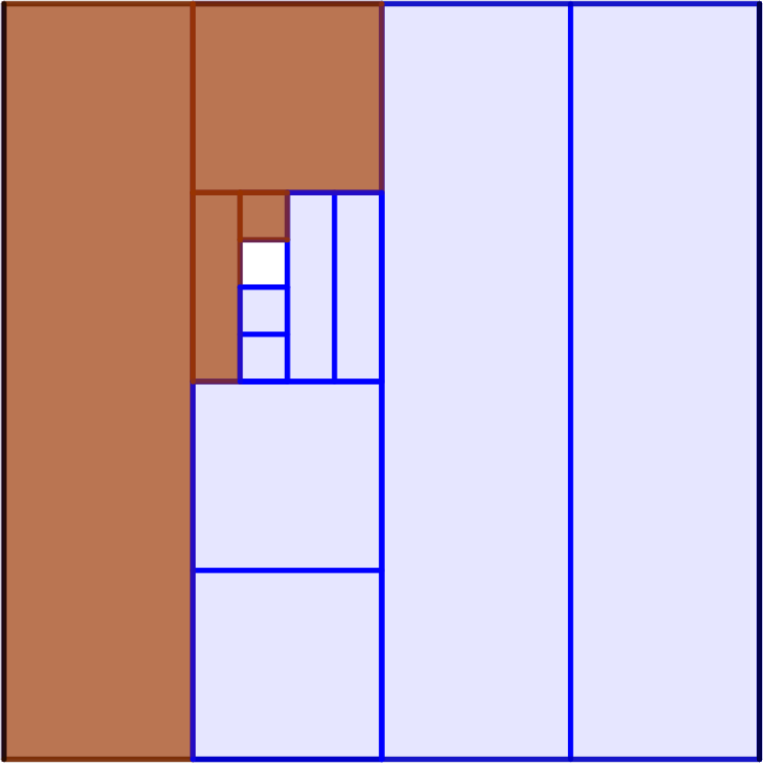}
\end{center}
\caption{A visual way of computing the geometric sum for $a=4$ for the iteration $n=3$
See also \url{https://www.geogebra.org/m/bt2vjyrv} for an interactive applet}
\label{chocolate}
\end{figure}

A mathematical introduction of number series was also started. A well-know concept of ``eating a table
of chocolate in infinitely many steps'' was used, similarly to Zeno's idea on ``Achilles and tortoise'' \cite{aristotle}.
\begin{enumerate}
\item You have a chocolate. Today you eat the half of it. Tomorrow you eat the
half of the rest. The day after tomorrow you eat again the half of the rest.
And so on. How much chocolate do you eat in 10 days? In 100 days?
\item You have a chocolate. Today you eat $1/2$. Tomorrow $1/4$. The day after
tomorrow $1/8$. And so on. How much chocolate do you eat in 10 days? In
100 days? What happens if you can eat the chocolate forever?
\item Same question, but you eat every day the one-third of the previous day.
\end{enumerate}
By using the hint that ``computing also the amount of chocolate you will never eat'' one can have a general
idea on computing the sum 
\begin{equation}\sum_{n=1}^\infty\frac{1}{a^n},
\end{equation} that is $1/(a-1)$ (see Fig.~\ref{chocolate}).
This idea can be useful to prepare the concept of self-similarity, hence the small square in the middle
is actually self-similar to the big one---if infinitely many iterations are considered.

\begin{figure}
\begin{center}
\includegraphics[width=0.6\textwidth]{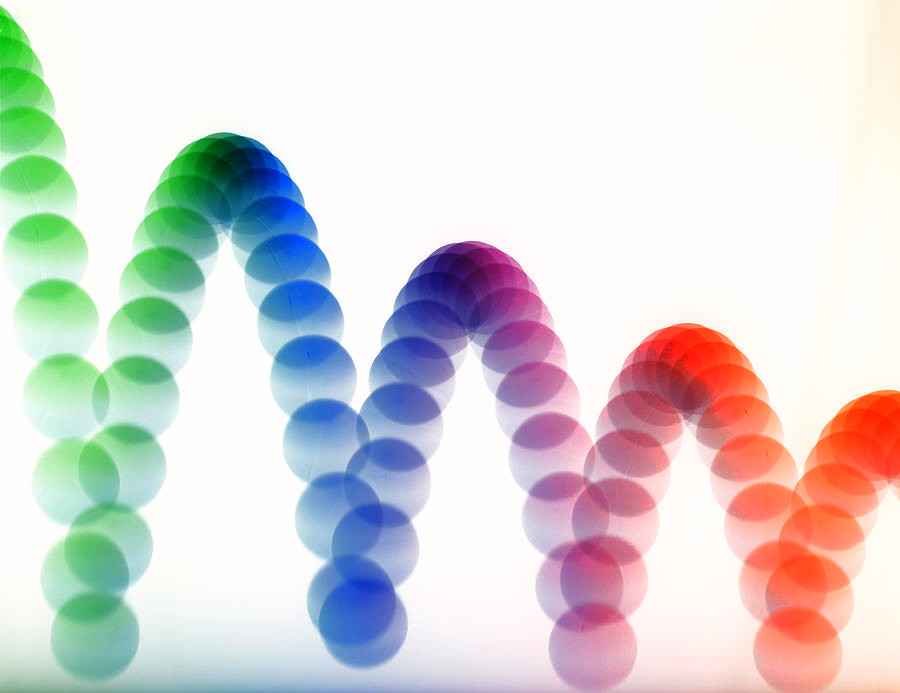}
\end{center}
\caption{Movement of a bouncing ball. This piece of artwork is based on a photo
provided by \url{https://stockcharts.com/articles/wyckoff/2015/08/follow-the-bouncing-ball-.html}}
\label{bouncing}
\end{figure}

\begin{figure}
\begin{center}
\includegraphics[width=0.6\textwidth]{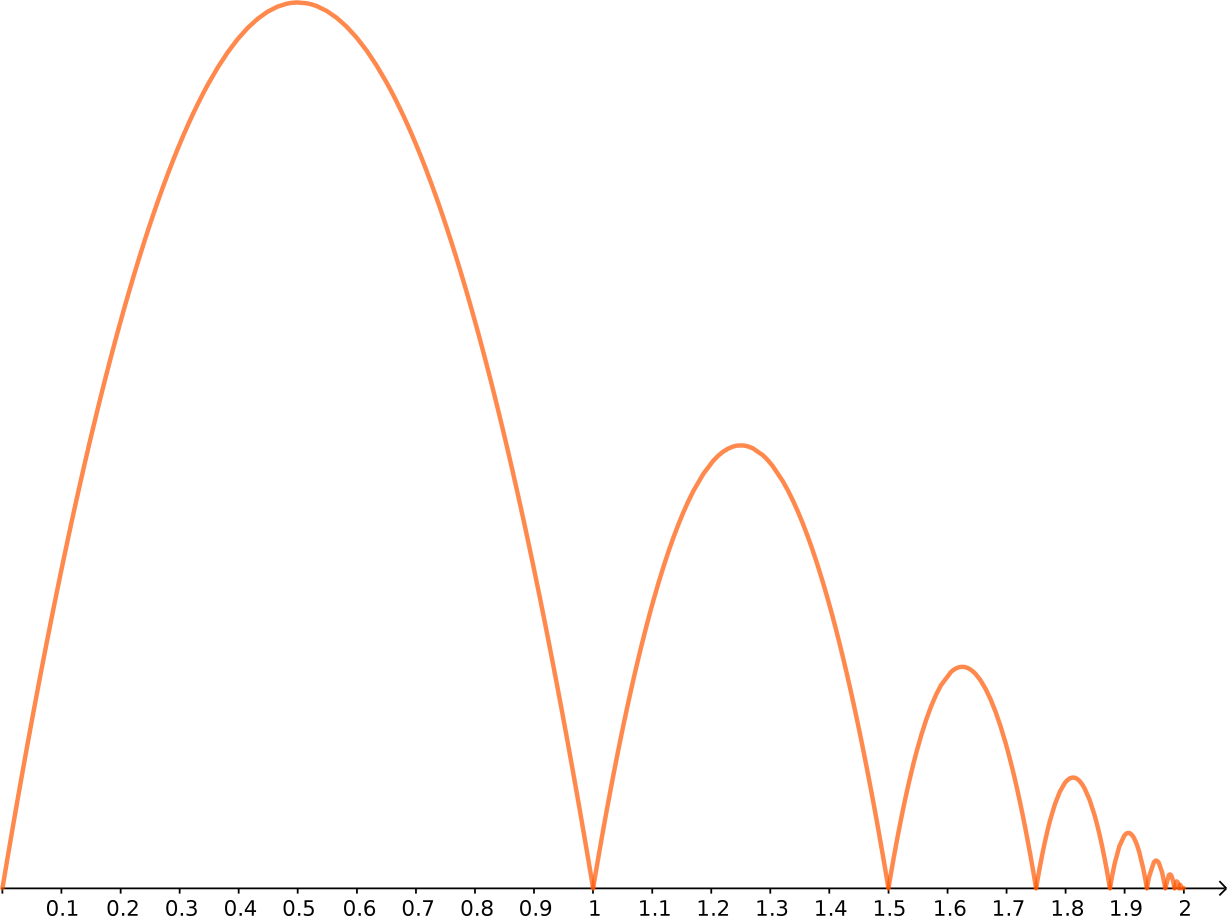}
\end{center}
\caption{Movement of a bouncing ball: an example of self-similarity}
\label{self-sim-par3}
\end{figure}

The movement of
a bouncing ball (Fig.~\ref{bouncing} and \ref{self-sim-par3}) can also be identified as an example of self-similarity
and the length of the movement can be computed by the geometric sum.

The final activity on the first day was to prove the sum of the geometric series algebraically.
Two proofs were discussed, one complete and an incomplete one. The incomplete one includes a step
that requires careful consideration of subtracting quantities that can be eventually infinite---this
fact was later discussed in detail. % Day ...

\begin{figure}
\begin{center}
\includegraphics[width=0.8\textwidth]{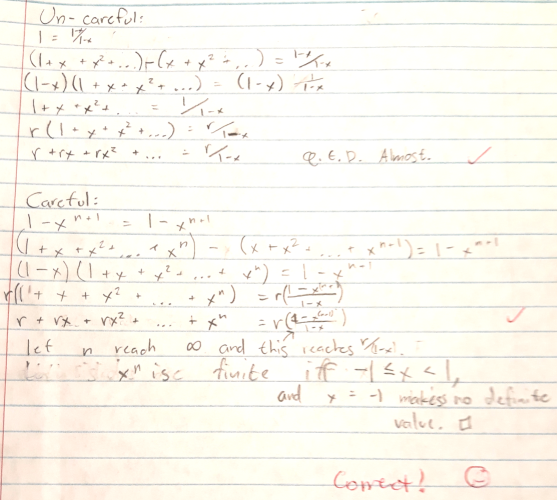}
\end{center}
\caption{A part of a student's worked-out homework for Day 1}
\label{hw1-3}
\end{figure}

One of the students already knew the formula of the geometric sum in the form
\begin{equation}\label{rmf}
\sum_{n=0}^\infty\frac{r}{x^n}=\frac{r}{1-x},
\end{equation}
and during the course we used this variant several times.
As a third homework activity the students had to give a detailed proof on (\ref{rmf}).
A worked-out homework can be seen in Fig.~\ref{hw1-3}. Seven students turned this homework
in, five of them were correct.

\subsection*{Day 2}

\begin{figure}
\begin{center}
\includegraphics[width=1\textwidth]{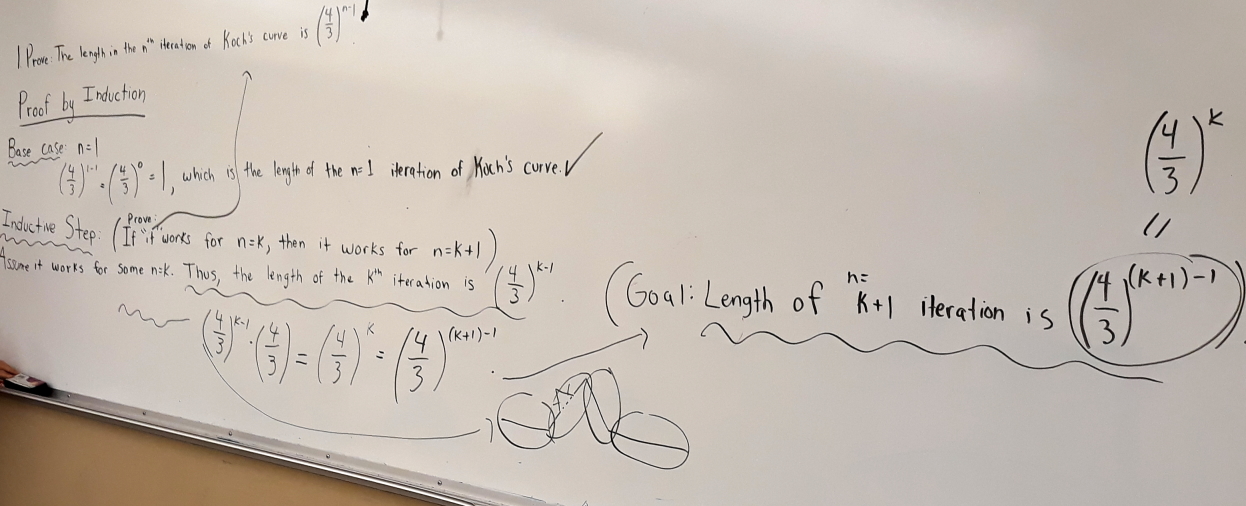}
\end{center}
\caption{An explained solution of a homework, made by the counselor Keri Celeste}
\label{Day2-Keri-edited}
\end{figure}

The Koch curve \cite{koch} was defined and the concept of iterations was introduced.
We discussed that the first iteration of the curve has length $1$, the second has $4/3$,
the third $4/3\cdot4/3$. By intuition we found the general formula $$(4/3)^{n-1}$$ for the
length of the $n$th iteration, but a rigorous proof---via induction---was left as homework.
(A rigorous proof is shown in Fig.~\ref{Day2-Keri-edited}.)

\begin{figure}
\begin{center}
\includegraphics[width=0.6\textwidth]{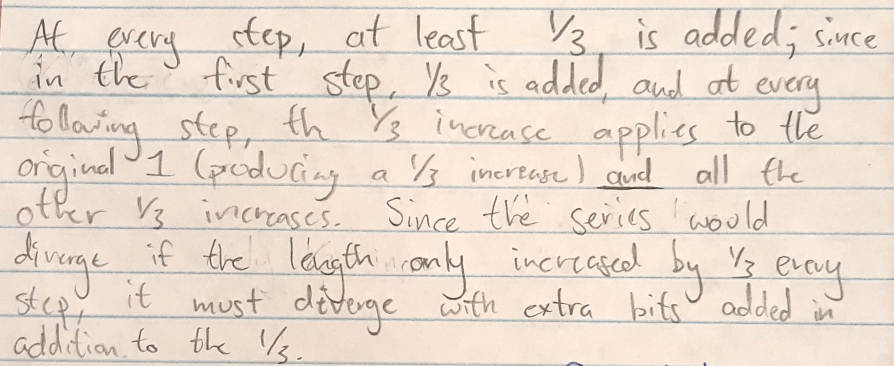}
\end{center}
\caption{A student's solution to prove that the length of the Koch curve diverges}
\label{Day2-Josh2-edited}
\end{figure}

It was mentioned that monotone increasing does not imply divergence to infinity,
but the length of the Koch curve actually diverges. A detailed proof on this was given as homework,
but first a similar question was discussed:
\begin{quotation}
There is a \$100 bond with a monthly interest rate 10\%, so one can gain \$10
after one month. In case of compound interest the bank pays not only \$10
after another month but \$11 because of the increased bond, which is \$110.
By using this concept, show that $100\cdot(1 + 1/10 \cdot m) \leq 100\cdot(1 + 1/10)^m$.
\end{quotation}
A student's solution on the divergence of the length of Koch's curve can be found in Fig.~\ref{Day2-Josh2-edited}.

\subsection*{Day 3}

\begin{figure}
\begin{center}
\includegraphics[width=0.4\textwidth]{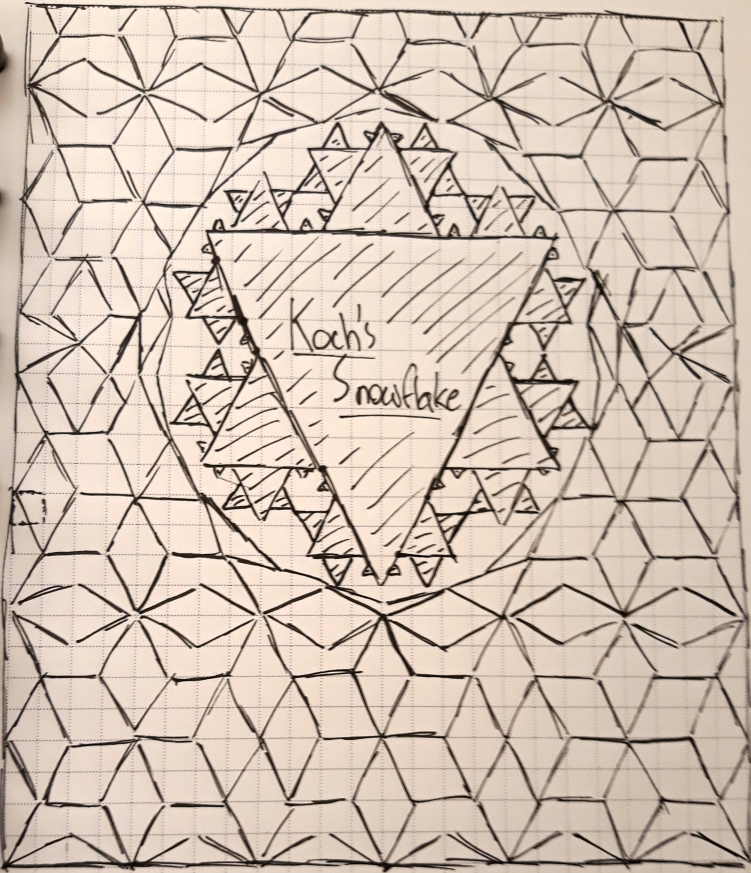}
\end{center}
\caption{Koch's snowflake as a student's artwork before Day 3}
\label{Day3-Max-snowflake-edited}
\end{figure}

\begin{figure}
\begin{center}
\includegraphics[width=0.4\textwidth]{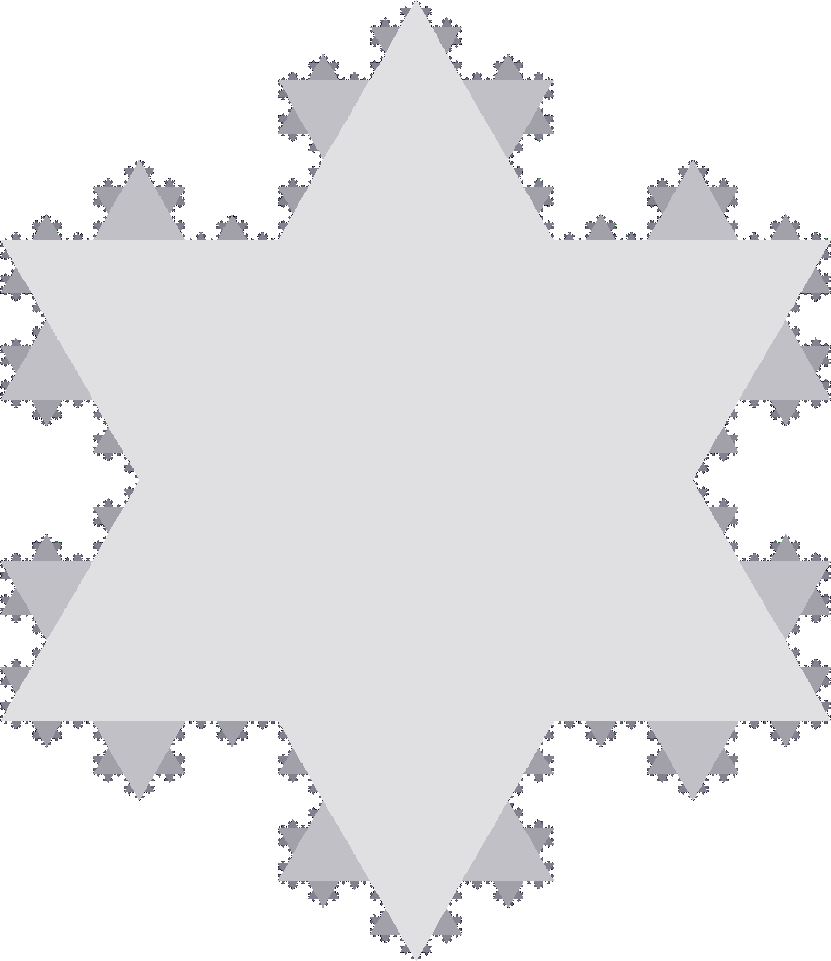}
\end{center}
\caption{Koch's snowflake shown in the software tool XaoS}
\label{koch-xaos}
\end{figure}

\begin{figure}
\begin{center}
\includegraphics[width=0.6\textwidth]{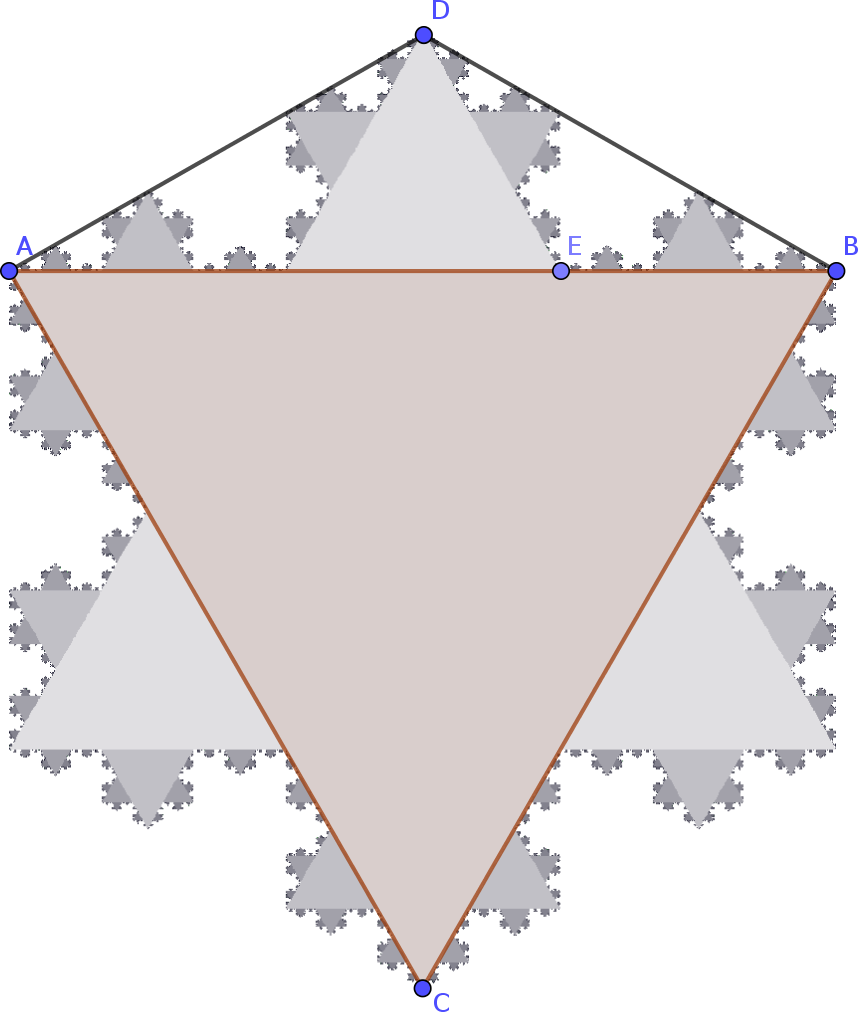}
\end{center}
\caption{Parts of Koch's snowflake are highlighted}
\label{koch}
\end{figure}

The software tool XaoS was introduced, limiting its use to display the Koch snowflake (Fig.~\ref{Day3-Max-snowflake-edited} and \ref{koch-xaos}).
We also computed its area by using (\ref{rmf}) after considering that the area of the internal
big triangle $ABC$ (see Fig.~\ref{koch}) is $T$, and then finding the sum
$$T+3\cdot\left(\frac13\right)^2\cdot T+
3\cdot\left(\left(\frac13\cdot\frac13\right)^2\cdot4\right)\cdot T+
3\cdot\left(\left(\frac13\cdot\frac13\cdot\frac13\right)^2\cdot4^2\right)\cdot T+\ldots
$$ that is
\begin{equation}
\begin{split}
T\cdot\left(1+3\cdot \sum_{n=0}^\infty{\left(\frac1{3^{n+2}}\right)^2\cdot4^n}\right)&=\\
=T\cdot\left(1+3\cdot \sum_{n=0}^\infty{\left(\frac1{3^2}\right)^{n+2}\cdot4^n}\right)&=\\
=T\cdot\left(1+3\cdot\frac19 \sum_{n=0}^\infty{\left(\frac49\right)^{n}}\right)&=\\
=T\cdot\left(1+\frac13 \cdot\frac{1}{1-4/9}\right)&=\\
=T\cdot\left(1+\frac13\cdot\frac95\right)&=\frac85\cdot T.
\end{split}
\label{T85}
\end{equation}

Without any computation, however, by intuition it is clear that the snowflake has a finite area, but this can also be directly
proven by induction when considering the following facts:
\begin{enumerate}
\item The first iteration of the snowflake is bounded (it is actually a triangle).
\item The second iteration of the snowflake is the first iteration extended by three triangles that each is bounded by
an isosceles triangle with angles $30^\circ$, $30^\circ$ and $120^\circ$. As a consequence, the second iteration is
bounded by a regular hexagon.
\item The third iteration is the second iteration extended by 12 triangles that each is bounded by an isosceles
triangle, similar to the one described above. For this reason all extensions are bounded by the isosceles triangles
given above, that is, the third iteration is still bounded by the above mentioned regular hexagon.
\item As a consequence, each iteration is bounded.
\end{enumerate}

The exact formulation of the proof was left as a challenge homework---that is not obligatory to solve.

Also, by using the idea from Fig.~\ref{chocolate}, another problem was to find a different kind of computation
of the area instead of the method described in (\ref{T85}).
To solve this, we denote the area of triangle $ABD$ by $t$, and that part of it which belongs to the snowflake we
call $U$. When assuming $T=1$ we search for the snowflake area $A=1+3U$. 
Now let us denote that area of the triangle $BDE$ that does not belong to the snowflake by $u$.
Clearly, $t=U+2u$. Geometrically it seems clear that the forms denoted by $u$ and $U$ are similar to each other,
and the similarity ratio is $1:3$, this follows from the fact that the area of triangle $BDE$ is one-third of
the area of triangle $ABD$. The final conclusion was left as homework.

\begin{figure}
\begin{center}
\includegraphics[width=0.6\textwidth]{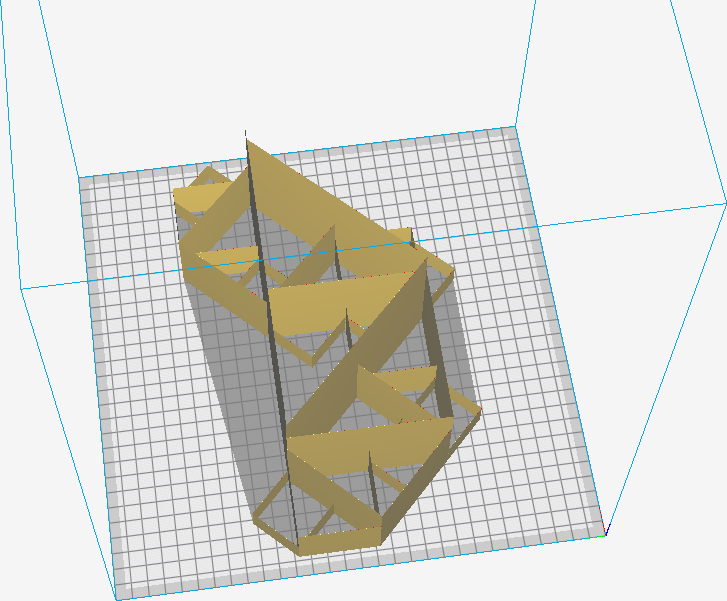}
\end{center}
\caption{A 3D object that contains more iterations of the dragon curve at the same time, shown in Ultimaker Cura
(see \url{https://ultimaker.com/software/ultimaker-cura})}
\label{cura-josh}
\end{figure}

One of the students gave a homework assignment to the author, namely, to find a way to print more than
one iteration of the dragon curve at the same time. This challenge was accepted, but to give an
answer the author needed a couple of days. (See Fig.~\ref{cura-josh} for an incomplete sketch.)

\subsection*{Day 4}
The completion of the homework was successful only for one student. So, we discussed the solution in the group as follows:
At the beginning we conclude that $t=3u+2u=5u$, and then, because for the regular triangle in the middle of $t$
the equation $t/3=T/9=1/9$ holds, and $t=1/3$, and therefore $u=1/15$, $U=1/5$, and finally $A=1+3\cdot1/5=8/5$ follows.

The main topic on Day 4 was to observe self-similarity on some other objects. An animation on pentagons
(see \url{http://prover-test.geogebra.org/~kovzol/pentagon/} and \url{http://prover-test.geogebra.org/~kovzol/pentagon/slider.html})
was shown. Simple questions like deciding self-similarity of a segment, a ray or a line were given
as classroom questions and homework.

\begin{figure}
\begin{center}
\includegraphics[width=0.6\textwidth]{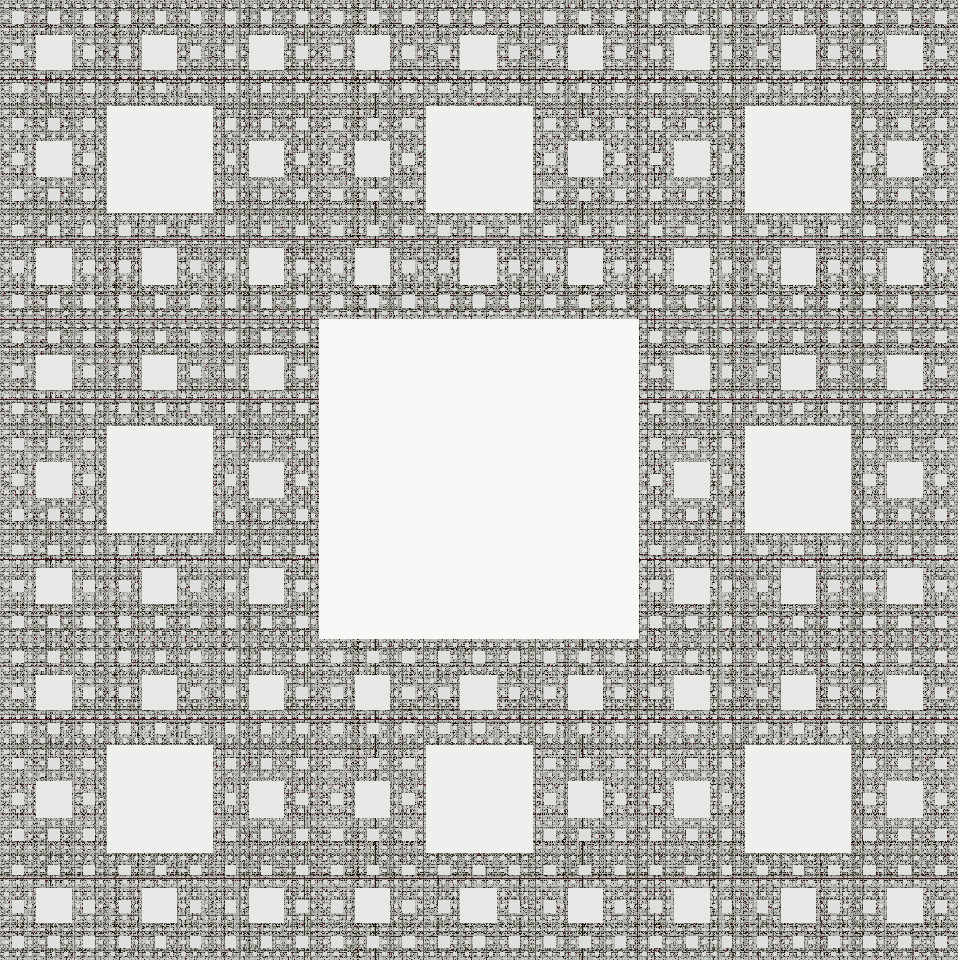}
\end{center}
\caption{Sierpi\'nski carpet in the software tool XaoS}
\label{scarpet-xaos}
\end{figure}

XaoS offers an easy way to study the Sierpi\'nski carpet \cite{sierpinski} (see Fig.~\ref{scarpet-xaos}, note that
some details of the picture are inaccurate due to numerical issues in binary-ternary conversions).
We computed that its perimeter is infinite but its area is zero---the latter fact was surprising for
many students. We used (\ref{rmf}) to prove the fact on the area, by computing
$$1-\left(\frac13\right)^2-8\cdot\left(\frac13\cdot\frac13\right)^2-8^2\cdot\left(\frac13\cdot\frac13\cdot\frac13\right)^2-\ldots$$

Students concluded the funny thing that creating a fence for the Sierpi\'nski carpet would cost an infinite amount of money,
but for this price the plot would not be worth anything.

To have more experience on proofs by induction, and to connect knowledge to ``textbook mathematics'', another
homework was to reformulate the ``monthly interest'' problem as Bernoulli's inequality in the form
$$(1+h)^n\geq1+hn$$ if $h\geq-1$, and induction had to be used for the proof.

Today's assignments were solved well by most of the students, but for the proof by induction the counselor's help was
required.

\subsection*{Day 5}

After Sierpi\'nski's carpet we observed its one-dimensional version, namely, the Cantor set \cite{cantor0,cantor}. We learned that similar
relations hold like in the two-dimensional case, namely, the length of the removed intervals is 1, therefore
Cantor's ternary set has a ``measure'' of zero. Here we did not introduce any exact definitions from set theory---the observed sums were
still analytically correct. A homework assignment was to finalize the computations precisely
by using (\ref{rmf}).

The other assignment was to express the numbers $1/2$ and $1/4$ in ternary form. Here again using (\ref{rmf}) was suggested.

Today's assignments were solved by about half of the students, but just a few solutions were completely correct.

Having the end of the week the web page \url{https://www.cutoutfoldup.com/216-dragon-curve.php} was advertised
to invite the students to try paper folding harder---and to be prepared for solving the puzzle on the dragon curve.

\subsection*{Day 6}

\begin{figure}
\begin{center}
\includegraphics[width=0.9\textwidth]{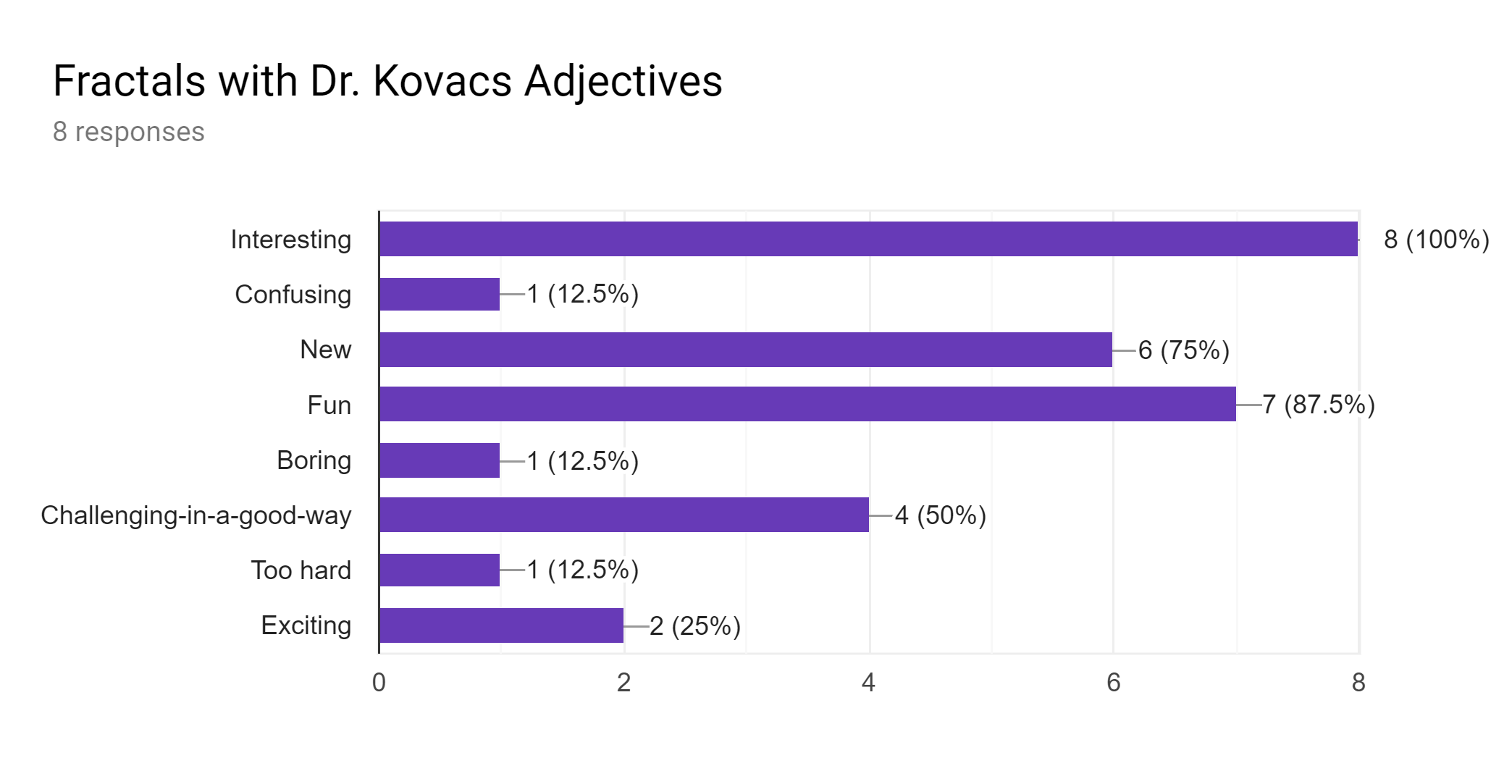}
\end{center}
\caption{Some responses of the students after the first week}
\label{adj1}
\end{figure}

At the end of the first week some feedback was given by the students (see Fig.~\ref{adj1}). 

We discussed Day 5's assignment on the ternary form of a number. To connect this concept to the Cantor set,
the following homework was given: \textit{Explain why the following statement is true: If the
ternary form of a number x contains the digit 1 then the number x is not a member of the Cantor set.}

Also, at the beginning of the first week a new topic was started: \textit{escape time fractals}.
We began the investigation with the number sequence 
\begin{equation}
(((x^2+x)^2+x)^2+x)^2+\ldots
\label{m}
\end{equation} for values $x=0,1,2,-1,-2$
and discussed the meaning of notions ``bounded'' and ``convergent''. A spreadsheet and a GeoGebra \cite{GG5} applet
\url{https://www.geogebra.org/m/r3c3vxcj} helped us visualizing
this `unusual' behavior. We experienced that for $0.26$ there is no convergence any longer, but for $0.25$ or
a bit below there is convergence. As a challenge homework assignment the students needed to show that the sequence (\ref{m})
is bounded if $x=0.25$ but it tends to infinity if $x>0.25$.

These assignments were not straightforward. However, one student solved both correctly.

\begin{figure}
\begin{center}
\includegraphics[width=0.45\textwidth]{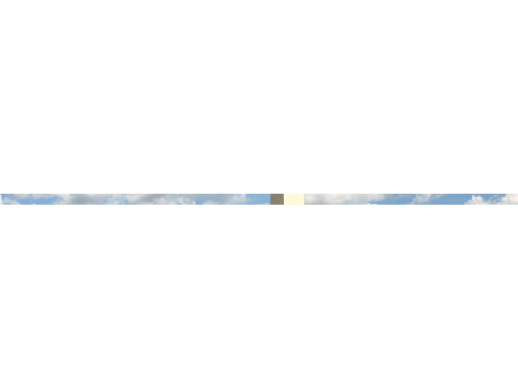}
\includegraphics[width=0.45\textwidth]{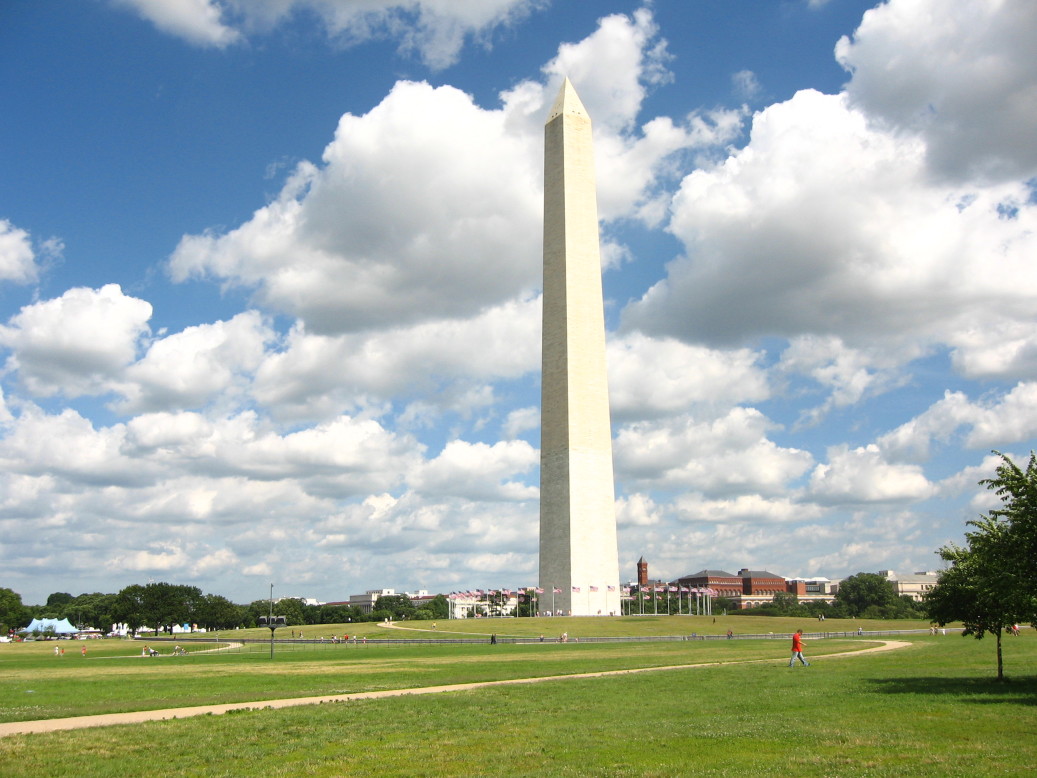}
\end{center}
\caption{A small stripe from the middle of a photo on the Washington Monument, and the full photo}
\label{wmon}
\end{figure}

As a final activity it was mentioned that a kind of generalization of this process will lead to the \textit{Mandelbrot set}
\cite{mandelbrot}.
The behavior constrained to the real numbers is something like knowing just one horizontal stripe of an image---we
cannot always tell the whole picture. (See Fig.~\ref{wmon} for an example.)

\begin{figure}
\begin{center}
\includegraphics[width=0.45\textwidth]{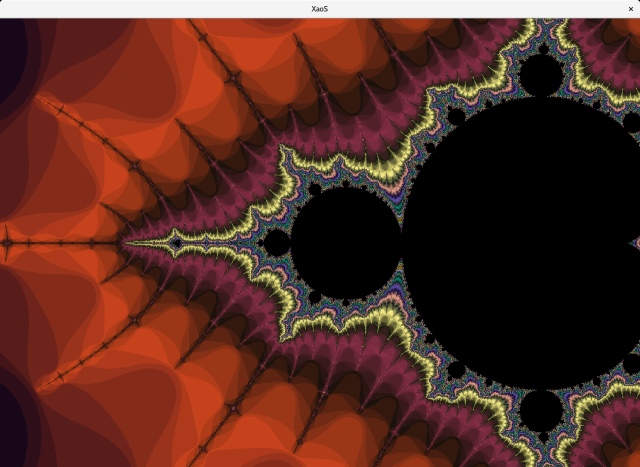} \includegraphics[width=0.45\textwidth]{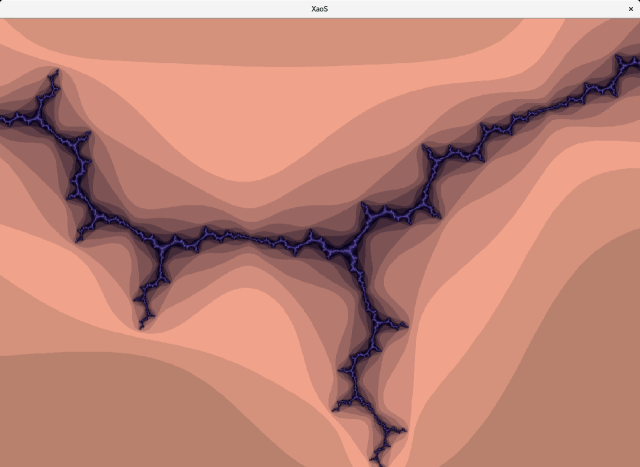}
\end{center}
\caption{Two screenshots of XaoS, generated by the students on Day 6}
\label{mbrot12}
\end{figure}

Following the students' request we finally played with XaoS to generate some beautiful fractal images on the Mandelbrot set
at the end of the class (see Fig.~\ref{mbrot12}).

\subsection*{Day 7}

First we discussed why (\ref{m}) is bounded if $0\leq x\leq1/4$. By setting 
\begin{equation}
x_1=x,\ x_{n+1}=x_n^2+x
\label{mseq}
\end{equation} we
can show that $x_n\leq1/2$. We use induction. For $n=1$ it is clear that $x_2\leq (1/4)^2+1/4=5/16<1/2$.
Now, if $x_n\leq1/2$, $x_{n+1}\leq(1/2)^2+1/4=1/2$.

On the other hand, another technique can be used to conclude that $x>1/4$ implies divergence.
If $x>0$, the sequence (\ref{mseq}) is monotone, and it is convergent if and only if for a certain number $c$
the equality
\begin{equation}
c^2+x=c
\label{q}
\end{equation}
holds. By solving the quadratic equation (\ref{q}) for $c$ we obtain
$$c=\frac{1\pm\sqrt{1-4x}}{2}$$
which is a real number if and only if $x\leq1/4$. This technique was mentioned on Day 6 quickly,
but in detail not given.

\begin{figure}
\begin{center}
\includegraphics[width=0.45\textwidth]{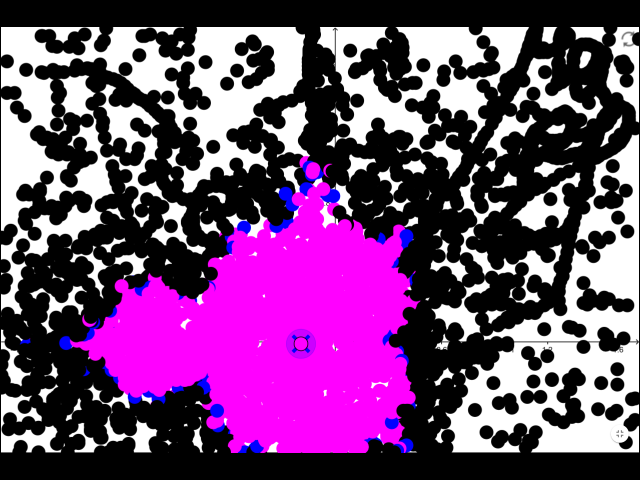}
\includegraphics[width=0.45\textwidth]{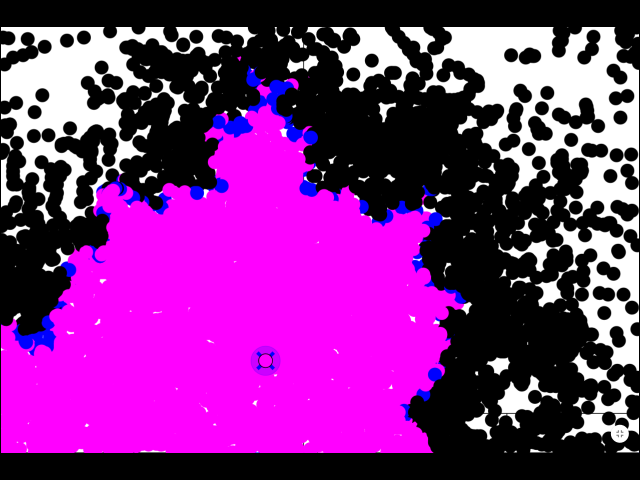}
\end{center}
\caption{Painting the Mandelbrot set}
\label{paint}
\end{figure}

The difficult theoretical details were made somewhat easier by fun activities. Students painted the Mandelbrot
set with the applet \url{https://www.geogebra.org/m/O2gLvtRY}. Fig.~\ref{paint} shows two screenshots
of such activities.

Instead of computing the sequence (\ref{mseq}) for certain complex numbers directly the applet
\url{https://www.geogebra.org/m/Npd3kBKn} was used. This helped to classify orbits that belong to
convergent, bounded or not bounded sequences. We explicitly defined the Mandelbrot set $\cal M$ by
considering all complex numbers $x$ such that (\ref{mseq}) is bounded.

We agreed on that the interval $[0,0.25]\subseteq\cal M$. We also checked that $-2\in\cal M$.
It seemed that for any positive $r$ the number $-2-r$ is no longer a member of $\cal M$. So we considered
the following sequence of expressions:
\begin{equation}
\begin{split}
&\;-2-r,\\
(-2-r)^2-2-r=&\;r^2+3r+2,\\
((-2-r)^2-2-r)^2-2-r=&\;r^{4} + 6 r^{3} + 13  r^{2} + 11  r + 2,\\
(((-2-r)^2-2-r)^2-2-r)^2-2-r=&\;r^{8} + 12 r^{7} + 62 r^{6} + 178 r^{5}+ 305 r^{4}\\
&\;  + 310 r^{3} + 173 r^{2} + 43 r + 2,
\end{split}
\label{43}
\end{equation}
and so on, the next element is
$r^{16} + 24  r^{15} + 268  r^{14} + 1844 r^{13} + 8726  r^{12} + 30012  r^{11} + 77290  r^{10} + 151258 
 r^{9} + 225873  r^{8} + 256068  r^{7} + 217186  r^{6} + 134202  r^{5} + 57809  r^{4} + 16118  r^{3} + 2541  r^{2} + 171  r + 2$.
%At this point we can recall some basic properties of convergent sequences, namely:
%\begin{theorem}
%\begin{enumerate}
%Given a number sequence $r_n\to0$. Then $r_n^n\to0$.
%Given two sequences $a_n\to A$ and $b_n\to B$. Then $a_n+b_n\to A+B$.
%Given a sequence $a_n\to A$ and a number $C$. Then $C\cdot a_n\to C\cdot A$.
%\end{enumerate}
%\label{t1}
%%\end{theorem}
%From Theorem \ref{t1} it immediately follows that if $r$ is ``close'' to $0$, then the sequence (\ref{mseq})
%stays in the neighborhood of $2$.

%On the other hand, to prove that $-2-r\not\in\cal M$ we do not need Theorem \ref{t1}.
By considering
the coefficients of $r$ in (\ref{43}) we find the following sequence:
\begin{equation}
3,11,43,171,\ldots
\label{s43}
\end{equation}

As homework assignment the students had to find a description for this sequence.
They could use an applet at \url{https://www.geogebra.org/m/yj5e9bng} to compute more
elements of the sequence.
As a challenge homework they had to show the divergence for $r>0$ as well.
Some parts of the first assignment was done properly by the half of the group---here some help was given them
in the office hours.

\begin{figure}
\begin{center}
\includegraphics[width=0.3\textwidth]{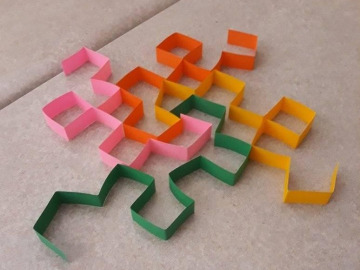} \includegraphics[width=0.3\textwidth]{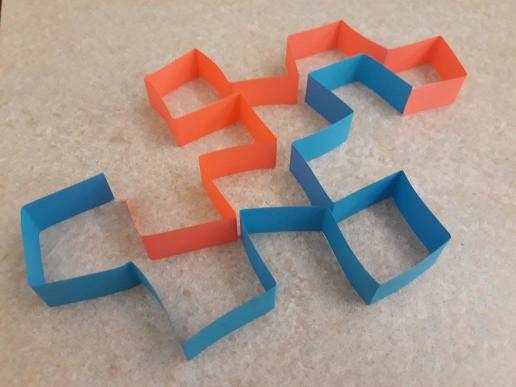}
\includegraphics[width=0.3\textwidth]{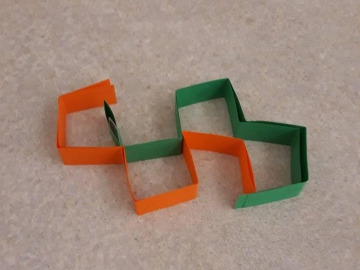} % !!! megvágni
\includegraphics[width=0.3\textwidth]{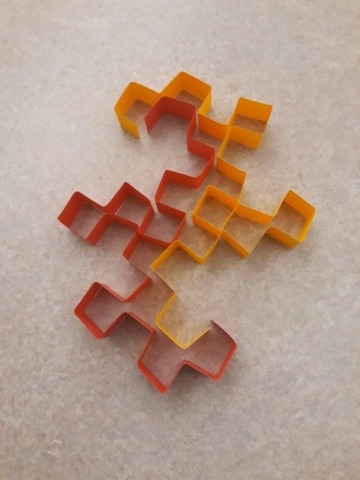} \includegraphics[width=0.3\textwidth]{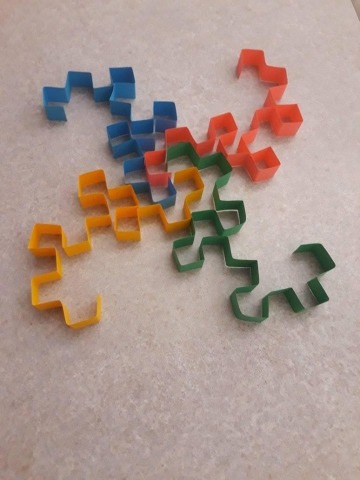}
\includegraphics[width=0.3\textwidth]{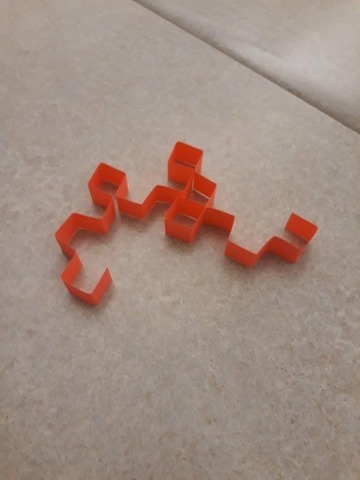} 
\end{center}
\caption{Dragon curves being fold out of paper}
\label{dragons}
\end{figure}

An evening activity dedicated to paper folding was advertised among the students.
One of them helped in folding a few first iterations of the dragon curve---the other students chose
a different evening activity. Later, however, all artifacts were exhibited in the students' classrooms,
so they were able to observe them by a closer look (Fig.~\ref{dragons}).

\subsection*{Day 8}

We discussed the solution of the ``normal'' homework first: A conjecture that $a_{n+1}=4\cdot a_n-1$
was found by several students, however, a proof was difficult to find. One of them, with some help, however,
produced the following derivation:

Assume that the $n$th expression is of form
$r^k+\ldots+a_n\cdot r+2$. Then, its square is of form
\begin{equation*}
\begin{split}
r^{2k}+\ldots+r^k\cdot a_n\cdot r\cdot + r^k\cdot2+ \\
+\ldots+\\
+ a_n\cdot r^{2k} + \ldots + a_n^2\cdot r^2+a_n\cdot r\cdot2+\\
+2r^k+\ldots + 2a_n\cdot r + 4.
\end{split}
\end{equation*}
Adding $-r-2$ and reading off the coefficient of $r$ we indeed get $2a_n+2a_n-1=4a_n-1$.

Now, we prove that $a_n\geq 3^n$. We use induction. For $n=1$ the statement is clear.
In case $a_{n+1}$ we have $a_{n+1}\geq 4\cdot3^n-1>3\cdot 3^n-1=3^{n+1}-1$. We conclude
that $a_{n+1}\geq3^{n+1}$.

During Day 8 we proved two more sophisticated theorems that required some basic knowledge of algebra
of complex numbers. In some sense this point was the most difficult part of the course since all
steps had to be described clearly, but a self-explanatory way was still preferred.

\begin{theorem}
If the sequence (\ref{mseq}) escapes from the origin further than two units, then it cannot go back.
\label{t1}
\end{theorem}

\begin{theorem}
If the sequence (\ref{mseq}) escapes from the origin further than two units, then it is not bounded.
\label{t2}
\end{theorem}

These theorems help defining a computer program to explicitly check if a complex number $x$ is
an element of the Mandelbrot set or not. We use the equality 
\begin{equation}
|u\cdot v|=|u|\cdot|v|
\label{uvprod}
\end{equation}
 that was
already known by a few students.
To highlight that we are working in the set of complex numbers, we denote the
start point by $z$ instead of $x$:
\begin{equation*}
x_1=z,\ x_{n+1}=x_n^2+z
\label{mseqz}
\end{equation*}

\begin{proof}
(Theorem \ref{t1}.) We consider two cases: \begin{enumerate}\item[(a)]First we assume that $|z|\leq2$. Let $E$ be a number
in the sequence (\ref{mseq}) such that $|E|>2$. That is $E$ is outside of the circle that has its center
in the origin and has radius 2. Because of the equality $|u\cdot v|=|u|\cdot|v|$ we learn that
$|E^2|=|E\cdot E|=|E|\cdot|E|=|E|^2>2\cdot2=4$, that is, $E^2$ is outside of the circle that has its center
in the origin and has radius 4. Now the next element in (\ref{mseq}) is $E^2+z$. Since $z$ is a vector
with a length less than equal 2, the point $E^2+z$ cannot be closer to the origin than 2 units---actually
it is further than 2 units.

\item[(b)]Second case: we assume that $|z|>2$. Consider the length of $z$, that is, $|z|=2+r$, $r>0$. Now let us compute
$|z^2+z|$. Since $|z|=2+r$, $|z^2|=(2+r)^2$, that is, $z^2$ lies on a circle that has its center
in the origin and has radius $(2+r)^2$. Now $|z^2+z|$ can be closer to the origin, but at most by
the length of $z$. We conclude that $|z^2+z|\geq|z^2|-|z|=(2+r)^2-(2+r)=r^2+3r+2$. Since $r^2+3r+2>2$,
the number $z^2+z$ is still outside the circle that has its center in the origin and has radius 2.
(Actually we implicitly used the triangle inequality in both cases.)
\end{enumerate}

(Theorem \ref{t2}.) We consider two cases again. First we assume that $|z|>2$, so let $|z|=2+r$. Now, by using
the previous thoughts we have the sequence $|z^2+z|\geq r^2+3r+2$,
$|(z^2+z)^2+z|\geq|(z^2+z)^2|-|z|$ (here we used the triangle inequality), and we get
$|(z^2+z)^2+z|\geq(r^2+3r+2)^2-2-r$. This recalls the same process as previously for the sequence $x=-2-r$
that diverges. By induction it is possible to see that the same formulas will appear. (We did not go into
the detail with the students, either.)

Second, we assume that $|z|\leq2$. Let $E$ be again a number
in the sequence (\ref{mseq}) such that $|E|>2$. Now $|E|=2+r$, $r>0$. Due to the triangle inequality $|E^2+z|\geq|E^2|-|z|
=(2+r)^2-|z|\geq(2+r)^2-2=r^2+4r+2.$ Here we either finalize the proof by pointing out that
this divergence is faster than for the first case since $r^2+4r+2>r^2+3r+2$, or simply use this last inequality
to give a proof based on the first case. As last step we highlight that the distance of $|E|$ from the
origin is at least $2+r$, the next element has a distance of $2+3r=2+r\cdot3$ at least, and by induction we can conclude
that the next distances will be at least $2+(r\cdot3)\cdot3$, $2+(r\cdot3\cdot3)\cdot3$, $\ldots$,
$2+r\cdot3^n\to\infty$.

\end{proof}

Here a graphical sketch was also created to illustrate the case (a) of the proof of Theorem \ref{t1}. The approach, however,
seemed to be unfamiliar and quite strange to the students. So proof of Theorem \ref{t2} was actually skipped and a new topic
was introduced, namely, \textit{Heron's process} (also known as Babylonian method, see \cite{fowler-rowson}),
as a preparation for the \textit{Newton fractal}\footnote{The first occurrence of the Newton fractal is difficult
to identify. Maybe \cite{peitgen-richter} was the first book that provided colored photos on various fractal images,
including the Newton fractal as well.}.

Today's homework assignments were:
\begin{enumerate}
\item Normal homework:
\begin{enumerate}
\item Consider the sequence $a_1=1$, $a_2=(a_1+2/a_1)/2$, $a_3=(a_2+2/a_2)/2$, $a_4=(a_3+2/a_3)/2$, $\ldots$\ Is
this convergent? Compute the first five elements to make a conjecture.
Using a calculator is allowed.
\item Prove (\ref{uvprod}) unless you already know how to prove it.
\end{enumerate}
\item Challenge homework:
\begin{enumerate}
\item We used the step $|z^2+z|\geq|z^2|-|z|$ in the proof of Theorem \ref{t1} (b) today. Explain why this step is
correct. Hint: Use the triangle inequality.
\item Prove Theorem \ref{t2}.
\end{enumerate}
\end{enumerate}

Assignment 1 (a) was solved by almost all students. About half of them managed to prove that the process converges
to the square root of 2. To achieve that, they used the technique from Day 7. The other assignments were not solved
and not even discussed during the remaining days.

\subsection*{Day 9}

\begin{figure}
\begin{center}
\includegraphics[width=0.7\textwidth]{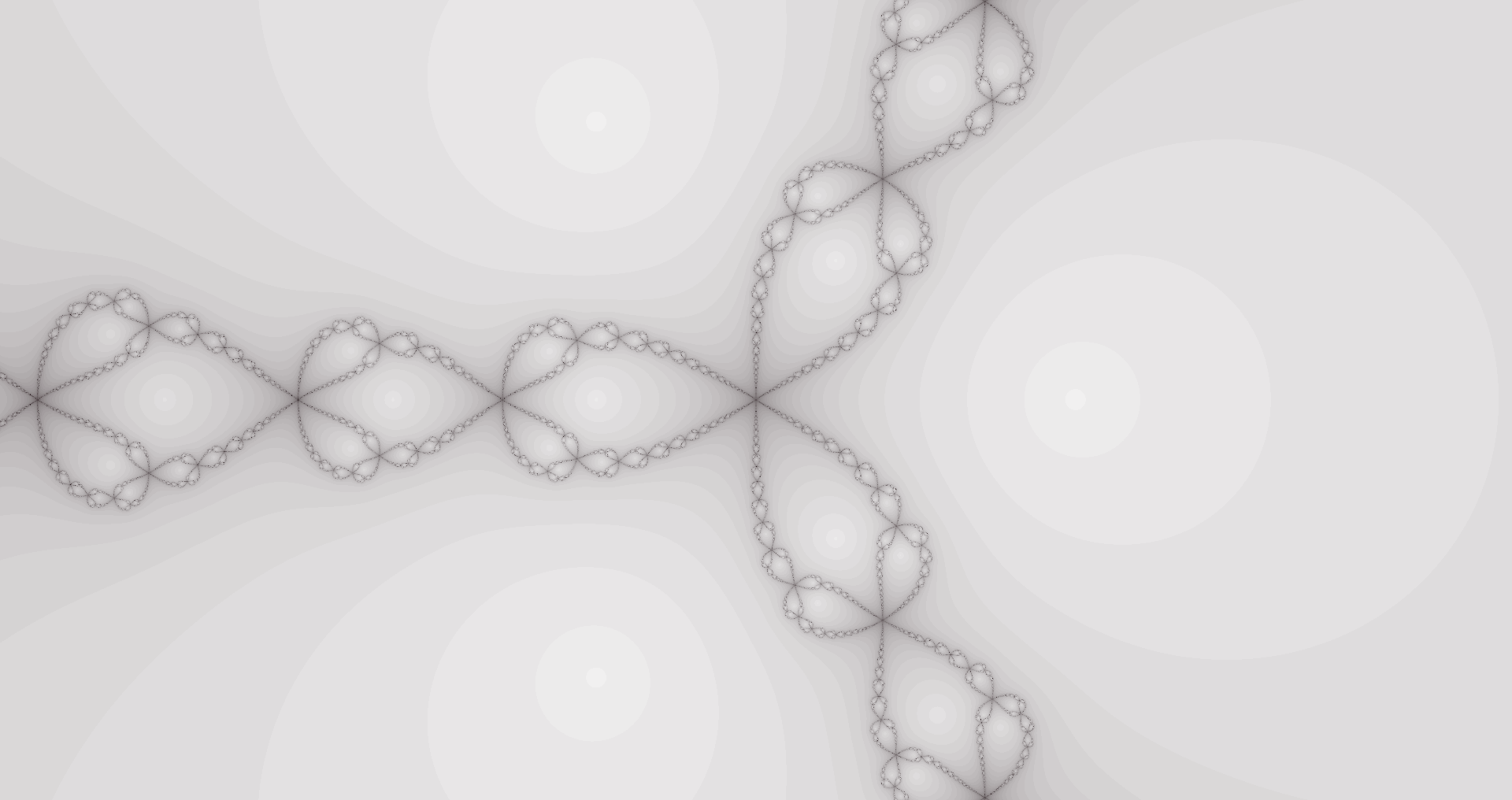}
\end{center}
\caption{Newton's cubic fractal in XaoS}
\label{newton-xaos}
\end{figure}

Today the work was dedicated to the introduction of Newton's cubic fractal (see Fig.~\ref{newton-xaos}). We generalized Heron's process
to compute the $k$th root of number $a$ approximately by using the sequence
\begin{equation}\label{newton}
x_1=z,\ x_{n+1}=A\left(\underbrace{x_n,x_n,\ldots,x_n,\frac{a}{x_n^{k-1}}}_{k\textrm{ numbers}}\right)
\end{equation}
where $z$ is a well-chosen number (for the first time it can be $1$, for instance),
$A(\ldots)$ is the arithmetic mean (that is, the average) of the numbers given in the parentheses.
We showed that if $z$ is a positive real number then the sequence is monotone after the second element.
An obvious lower bound of the sequence is $0$, but we showed that $x_n\to \sqrt[k]{a}$ is also a lower
bound, by using the AM-GM inequality for $k$ numbers.
Also, we concluded, by using the theorem on a monotone bounded sequence that $x_n\to \sqrt[k]{a}$.

Some parts of this journey were given as homework assignments. In particular, a Wikipedia article
\url{https://en.wikipedia.org/wiki/Inequality_of_arithmetic_and_geometric_means#Proof_by_Cauchy_using_forward%E2%80%93backward_induction}
was suggested as good read
that explains Cauchy's proof for the AM-GM inequality in general---for $2$ numbers this fact was
already known for most students, including the proof as well. 

During the journey a student came up with the idea to use the formula
\begin{equation}
x_{n+1}=A\left(x_n,\frac{a}{x_n^2}\right)
\end{equation}
to compute $\sqrt[3]{a}$. By using a spreadsheet we experienced that his suggestion seems correct.
However, we did not prove or disprove this idea.
%=($B6+$B$1/($B6^2))/2

The main concept of Newton-type fractals is to find all numbers $z$ in (\ref{newton}) such that
the sequence does not converge. This happens, of course, if there is a division by zero in the process.
Otherwise the process should be successful and it can deliver the $k$th root of $a$ approximately
after a few steps---here we did not discuss the speed of convergence. For real
numbers it is usually clear which number is the $k$th root of $a$. To generalize this concept we
observed the case $k=3$, $a=-1$. As a normal homework assignment, students had to verify that
$$\frac12-\frac{\sqrt3}2\cdot i$$ is a cubic root of $-1$ (the other complex root was checked
in the classroom already). A challenge homework was to find the non-real cubic roots of $-4$.
 
The normal homework assignments were solved by the majority of students.

\subsection*{Day 10}

During the last day two topics were covered: first, the explanation of the structure of the
cubic Newton fractal, and second, a sketch of a proof of two properties of the dragon curve
was shown.

\begin{enumerate}

\item A geometrical meaning of extracting the $n$th root of a complex number was explained.
For this, first the trigonometrical form of a complex number was introduced, but without
using $\cos$ and $\sin$ directly. Instead, the ``angle of a planar vector'' was used,
and it was shown, that for two complex numbers $u=a+b\cdot i$, $v=c+d\cdot i$, where both
lengths are $1$, the angle of the product of them is the sum of their angles. That is,
\begin{equation}
\arg u\cdot v\equiv\arg u+\arg v
\label{arg}
\end{equation}
where $\arg$ stands for the ``argument'' (in other words, ``angle'') of a complex number
and $\equiv$ allows angles like $0^\circ$ and $360^\circ$ to be considered equal.
(See the GeoGebra applet \url{https://www.geogebra.org/m/utqSs35J} for a sketch that actually explains
the trigonometric formulas $\cos(\alpha+\beta)=\cos\alpha\cos\beta-\sin\alpha\sin\beta$
and $\sin(\alpha+\beta)=\sin\alpha\cos\beta+\cos\alpha\sin\beta$---the applet is based on \cite[p.~46]{nelsen2})
Clearly, by using (\ref{uvprod}) as well, it is possible to study root extraction of a number
that lies on the unit circle by considering only rotations. In particular,
we determined all three complex roots of $i$ by finding all possible points of the unit circle
that satisfy the equation $\arg u^3=\arg i$. By using (\ref{arg}) we obtain that $3\arg u\equiv\arg i$.
In fact, the equation (\ref{arg}) in the given form was not shown to the students, but just
a verbal description was provided due to lack of time. It seemed however straightforward
that a complex number having unit length has three cubic roots.

But the same fact is true for any complex number, unless it is zero. In other words, the equation
$$z^3=a$$ for a given $a$ has always three solutions (unless $a=0$). For the final considerations
of this topic we will need the general form of this fact, namely,
\begin{theorem}
Any cubic polynomial has three complex roots (with multiplicities).
\label{fta3}
\end{theorem}
Or, an even more general form of this is the Fundamental Theorem of Algebra\footnote{Its
first rigorous proof was published in \cite{cauchy}, but the proof was actually found by Argand in 1806}, namely,
\begin{theorem}
Any polynomial of degree $k$ ($k>0$) has $k$ complex roots (with multiplicities).
\label{fta}
\end{theorem}
During the class we gave a sketch of a topological proof of Theorem \ref{fta}
by using the GeoGebra applet at the page \url{https://www.geogebra.org/m/HV6iUi3C}. The students
were pointed to learn more on this topic at \url{https://en.wikipedia.org/wiki/Fundamental_theorem_of_algebra#Topological_proofs},
in particular in section ``another topological proof''.

At this point all ingredients were ready to describe the structure of the cubic Newton
fractal, and, in fact, for all Newton-type fractals for higher degrees as well. We considered,
however, only the cubic case. The formula (\ref{newton}) can always be evaluated unless
$x_{n}=0$. For simplicity, let us assume $k=3$ and $a=8$. Also, assume that the next iteration fails,
that is $x_{n+1}=0$. Then, the iteration in (\ref{newton}) can be written as
\begin{equation*}
%\begin{split}
0=x_{n+1}=2x_n+\frac{8}{x_n^2}
%\end{split}
\end{equation*}
that is equivalent to
$$2x_n^3+8=0$$ that implies $$x_n=\sqrt[3]{-4}$$
that has one real root: $-\sqrt[3]4$, and two complex roots:
$\sqrt[3]4\cdot\left(1\pm\frac{\sqrt3}2\cdot i\right)$. This identifies three knots
in the cubic Newton fractal. That is, if $x_1$ is set to one of these three numbers,
$x_2$ will be zero and the iteration fails. Clearly, $x_1=0$ is also a bad choice
(in this case the iteration process stops immediately), so $0$ is actually the very
first identified knot.

But assume now that the number $c$ is an already identified knot, that is, it must
be avoided, unless in some steps the iteration runs into division by zero.
Then (\ref{newton}) can be
written as
\begin{equation*}
%\begin{split}
c=2x_n+\frac{8}{x_n^2}
%\end{split}
\end{equation*}
that is equivalent to
$$c\cdot x_n^2=2x_n^3+8,$$ that is
$$2x_n^3-c\cdot x_n^2+8=0$$ which is a cubic equation that has three complex roots
$c_1$, $c_2$ and $c_3$
(with multiplicities), according to Theorem \ref{fta3}. This simply means that
the knot $c$ introduces three additional knots $c_1$, $c_2$ and $c_3$. In fact,
the knot $0$ introduces three knots (the complex cubic roots of $-4$), and each
introduces three more knots. This process repeats itself in an infinite pattern.

The first knots were computed with GeoGebra as well (Fig.~\ref{badguys}). The relationship
between the knots are visualized by using an artifact in Fig.~\ref{trif} that was
partially reconstructed as a Blender image (Fig.~\ref{thething}) and then exported as an STL file, and finally
printed by a 3D printer. Lastly, the knot structure is shown if Fig.~\ref{dot}.

\begin{figure}
\begin{center}
\includegraphics[width=0.8\textwidth]{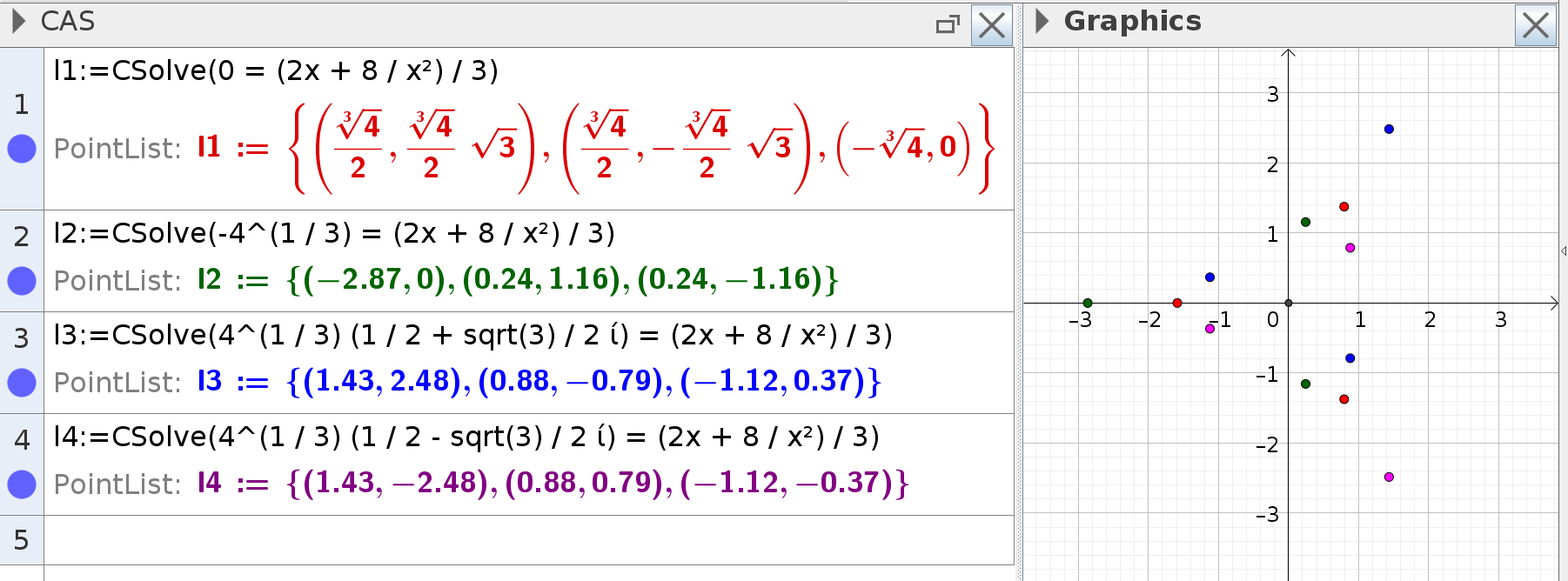} 
\end{center}
\caption{The first knots of the cubic Newton fractal: $0$ is black, the cubic roots of $-4$
are red, and the additional knots they introduce are green, blue and magenta, respectively}
\label{badguys}
\end{figure}

\begin{figure}
\begin{center}
\includegraphics[width=0.5\textwidth]{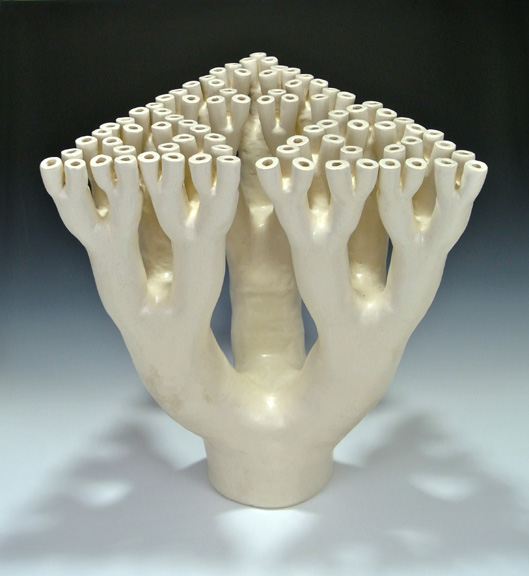} 
\end{center}
\caption{``Trifurcation'', a sculpture made by Robert Fathauer in 2014, \url{http://robertfathauer.com/Trifurcationb.html}.
Note that each iteration has a layout like an iteration of the Sierpi\'nski triangle}
\label{trif}
\end{figure}

\begin{figure}
\begin{center}
\includegraphics[width=0.5\textwidth]{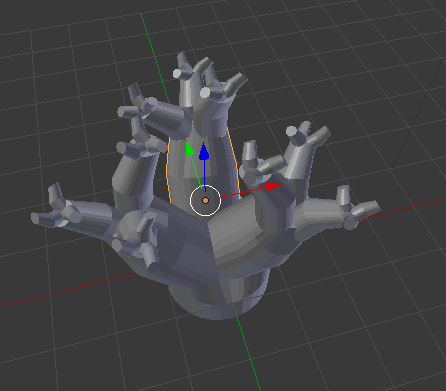} 
\end{center}
\caption{Benedek Kov\'acs's ``The Thing'', a partial reconstruction of Fathauer's sculpture in Blender}
\label{thething}
\end{figure}

\begin{figure}
\begin{center}
\includegraphics[width=\textwidth]{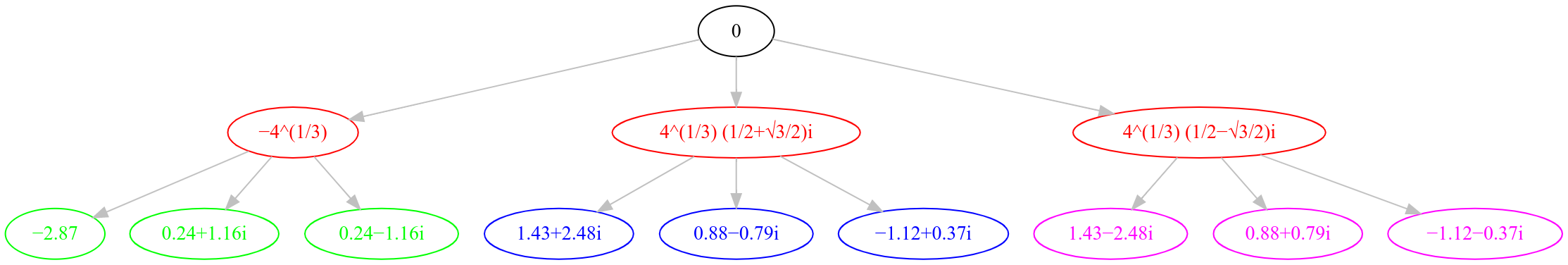} 
\end{center}
\caption{The beginning of the knot structure of the cubic Newton fractal}
\label{dot}
\end{figure}

%http://robertfathauer.com/Trifurcation.html

\item The second topic was to finalize discussing some facts on the dragon curve that were mentioned during
the first days. The dragon curve has several interesting properties, but at camp we focused on two statements:
\begin{enumerate}
\item The dragon curve is not self-intersecting and not self-overlapping.
\item Four copies of the dragon curve tile the plane without gaps.
\end{enumerate}
The second statement needs some explanation. Here we mean that if the plane is defined as an arbitrary large
(but finite) square, there is a certain iteration of the dragon curve that, considering 4 copies of it, will
cover the whole square (practically all edges in the rectangular grid).

Both questions are studied in detail in the literature, and during the course, we just mentioned the main
idea how one can prove such statements.

An excellent and exhaustive explanation of several features of the dragon curve is given in \cite{ryde-dragon}
by Kevin Ryde. Here we
actually just point the reader to consult this reference in Section 1.1 (Plane Filling) that proves
both theorems (both were communicated by Davis and Knuth first in \cite{knuth-davis1} and \cite{knuth-davis2}):

\begin{itemize}
\item The dragon curve touches at vertices but does
not cross itself and does not overlap itself.
\item Four copies of the dragon curve arranged at
right angles fill the plane.
\end{itemize}

Ryde gives a couple of ways to communicate the main idea of the proof. During the course
one of those communications, a classical tiling pattern K02A (see \url{http://tilingsearch.org/HTML/data24/K02A.html})
was used that is originated to Edgar (see \cite{edgar}). There was a convenient way to use this kind of
tiling: the campus of University of Colorado, Colorado Springs, has several outdoor public areas that have the same
type of tiling on the floor (see Fig.~\ref{grid1}). A simple way to explain the concept is to use these tiles actually in front of
the building of the lecture hall, and to draw on the tiles with various colored chalks. Unfortunately, this
plan was not achieved because of bad weather and lack of colored chalks---finally the explanation was performed
in the classroom by projecting the patterns on the whiteboard and extending them by colored pens.

\begin{figure}
\begin{center}
\includegraphics[width=0.4\textwidth]{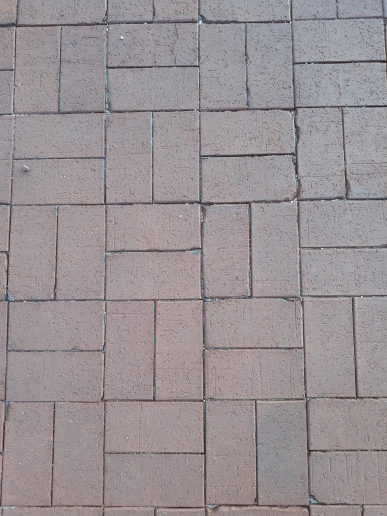} 
\end{center}
\caption{A K02A type grid at University of Colorado, Colorado Springs}
\label{grid1}
\end{figure}

By following \cite{ryde-dragon} and \cite{edgar} we consider a scaled and rotated copy of Fig.~\ref{grid2}
and observe it and the original picture at the same time. The ratio of scaling is $\sqrt2:1$ and the rotation
angle is $45^\circ$ clockwise.

\begin{figure}
\begin{center}
\includegraphics[width=0.4\textwidth]{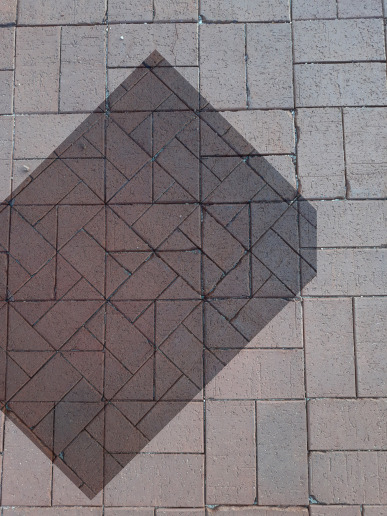} 
\end{center}
\caption{Considering Fig.~\ref{grid1} and its scaled-rotated variant at the same time}
\label{grid2}
\end{figure}

Let us consider now Fig.~\ref{grid4}. The second iteration of the dragon curve is drawn in red
on some edges of the first tiling, namely on those that belong to the edge of a square.
Similarly, the third iteration is shown in blue on the scaled-rotated tiling.

\begin{figure}
\begin{center}
\includegraphics[width=0.4\textwidth]{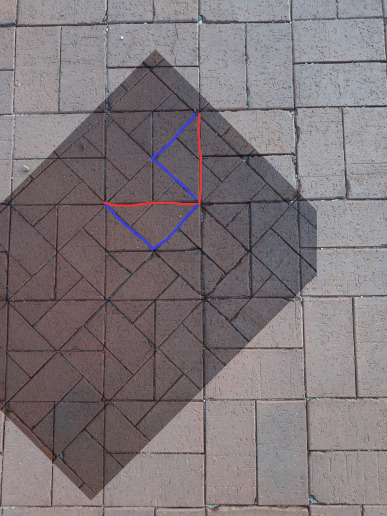} 
\end{center}
\caption{The second and third iterations of the dragon curve in red and blue, respectively}
\label{grid4}
\end{figure}

A general rule in this way of iteration is that a segment will be expanded into two segments and
they occupy that rectangle of the first grid which is ``free''. The blue segments are drawn
as edges of the second grid. After another iteration the blue segments
expand into green ones by following the same rule (see Fig.~\ref{grid5}): they occupy that
rectangle of the second grid which is ``free''. On the other hand, the green segments will be
drawn as edges of a third grid which is again a scaled-rotated variant of the second grid:
by magnification factor $\sqrt2:1$ and rotation angle $45^\circ$ clockwise (not shown in the figure).

\begin{figure}
\begin{center}
\includegraphics[width=0.4\textwidth]{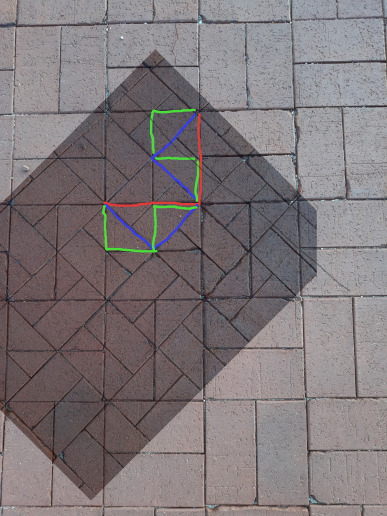} 
\end{center}
\caption{The second, third and fourth iterations of the dragon curve in red, blue and green, respectively}
\label{grid5}
\end{figure}

To confirm that the dragon curve does
not cross and overlap itself we can observe two other beginnings of iterations in Fig.~\ref{grid9}, one
in yellow and one in cyan. These examples are squares. Both expand to a union of two squares, but in a slightly different way,
into orange and darkcyan colors.

\begin{figure}
\begin{center}
\includegraphics[width=0.4\textwidth]{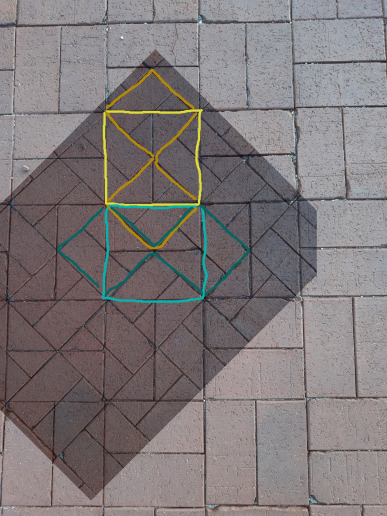} 
\end{center}
\caption{The second and third iterations of two variants of squares}
\label{grid9}
\end{figure}

We note that the first iterations overlap each other only at their joint edge. The second iterations
overlap each other only at the images of the overlapping edge of the first iteration, and have no self-intersections.
Now, by considering
the whole grid as another example of a first iteration, its second iteration will be still non-overlapping
and non-intersecting because of the same reason as mentioned for Fig.~\ref{grid9}. This property remains during each iteration.
As a consequence, any iteration of the whole grid remains non-overlapping and non-intersecting.

Finally, when observing the dragon curve we learn that its first iteration is a subset of the whole grid.
But the non-intersecting and non-overlapping properties hold for any subsets as well. Therefore, any iteration
of the dragon curve is non-intersecting and non-overlapping.

Before continuing with the second property, we highlight that the two first squares in Fig.~\ref{grid9}
initiate a kind of tiling of the grid. By observing two more squares, that is, a start with $2\times2$ squares
(as seen on the left part of Fig.~\ref{tile48} in red),
they will be expanded into a pattern (as seen in blue) that tiles the plane (as seen in Fig.~\ref{tiling}).
Clearly,
this tiling is non-intersecting and non-overlapping. 
%On the other hand, by considering the second iterations
%of the two first squares we learn that those 2+2 squares initiate another tiling of the grid, and that is
%again non-intersecting and non-overlapping---because .
Actually we used this property during the proof, and also the idea, that the tiling pattern contains nothing but just
a set of squares.

\begin{center}
*
\end{center}

The second property uses a simpler idea, but it can be visualized on the K02A tiling again. Fig.~\ref{tile48}
shows the first, second and third iterations of a $2\times2$ grid
in red, blue and green, respectively. One can learn that the third (green) iteration almost fills
a square that has 4 unit long sides---here we mean a unit according to the same iteration. During
the course we did not verify (and we leave to the reader here, too) that the 4th iteration fills
a $4\times4$ square completely, but the immediate consequence of this fact is that the 7th iteration
fills a $8\times8$ square completely, and so on.

Since the first iteration can be reached in two iteration steps from a start with four copies of the first
iteration of the dragon curve---see Fig.~\ref{dragons} for a quick verification---it is clear that
the later iterations fill an arbitrary big square.

Both properties were proven during the class only by sketching up the main ideas. The rest---to fill the gaps in---%
remained as informal homework assignments.

\begin{figure}
\begin{center}
\includegraphics[width=0.4\textwidth]{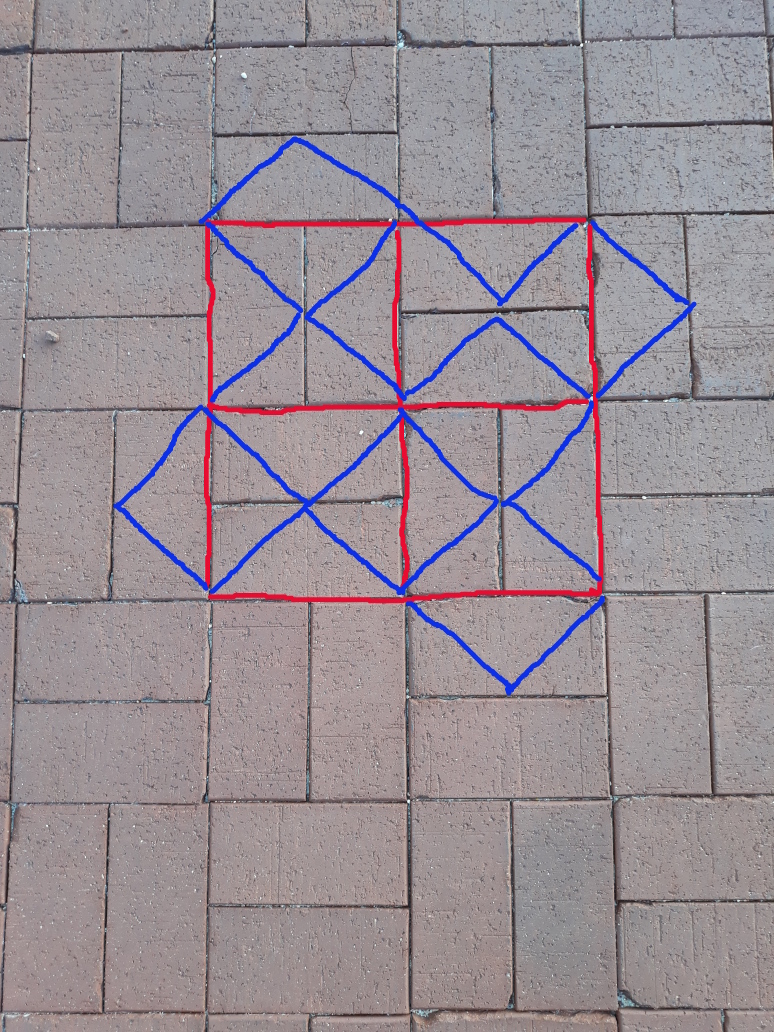} 
\includegraphics[width=0.4\textwidth]{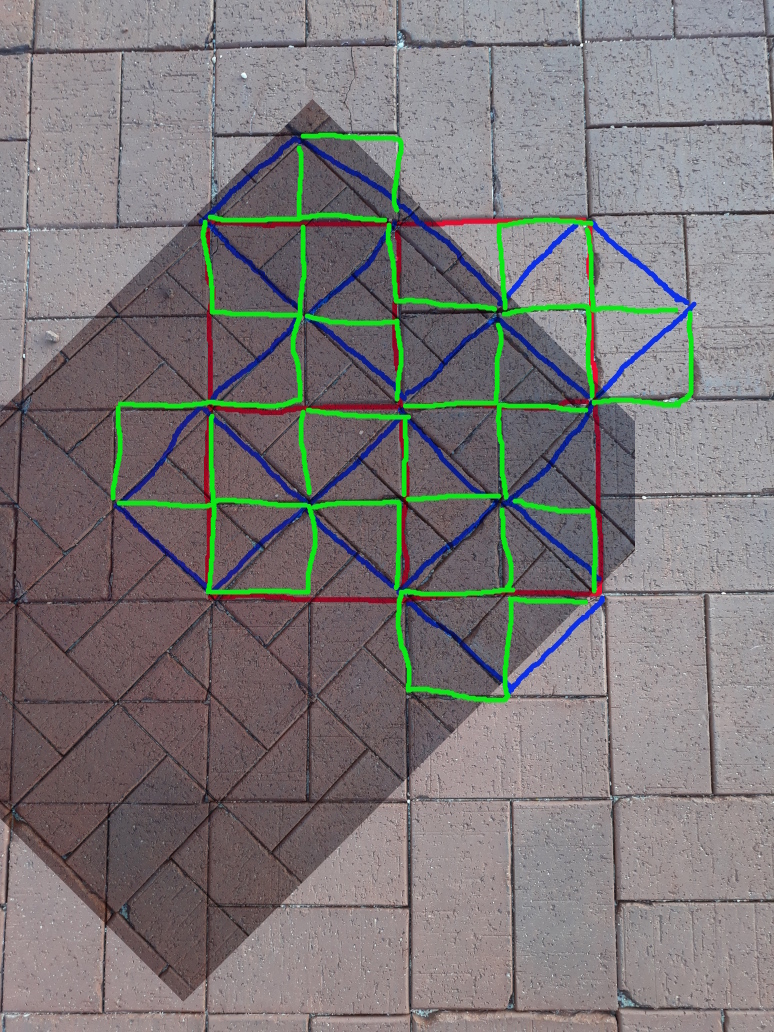} 
\end{center}
\caption{The 1st, 2nd and 3rd iterations of a $2\times2$ square}
\label{tile48}
\end{figure}

\begin{figure}
\begin{center}
\includegraphics[width=0.6\textwidth]{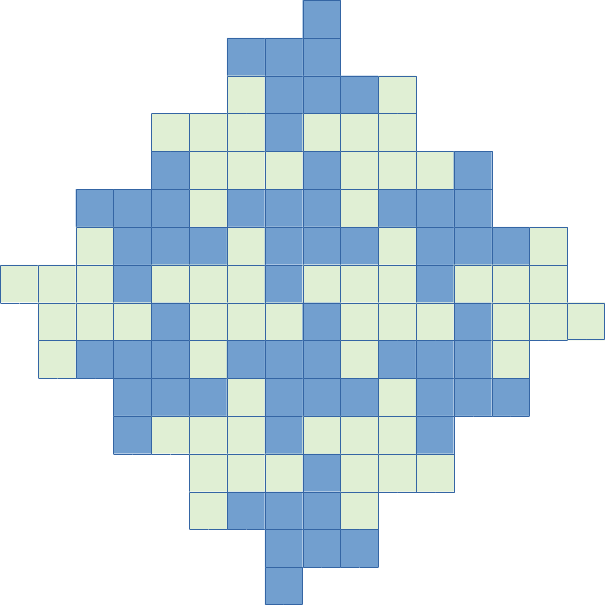} 
\end{center}
\caption{Starting the iteration with a bigger square the $2\times2$ parts will be expanded to a specific pattern}
\label{tiling}
\end{figure}

\end{enumerate}

%I took the proof sketches about the non-intersecting and non-overlapping
%property and the plane-tiling of the dragon curve from
%https://download.tuxfamily.org/user42/dragon/dragon.pdf. This is a very
%good read, especially the first two theorems.

%The proofs are illustrated with a tiling that is also visible at the UC
%campus. I uploaded all illustrations to
%http://test.geogebra.org/~kovzol/talks/epsiloncamp/Fractals/ (folders
%"grid" and "tile").

%There were some technical issues yesterday when plotting the graph of the
%"bad guys" of the cubic Newton fractal in GeoGebra. Here is a fixed
%version attached. Computing the cubic roots does not work in a web
%browser currently so one needs GeoGebra 5 or 6 Classic installed locally
%to load the .ggb file.

%badguys.ggb

\section{Conclusion}

\begin{figure}
\begin{center}
\includegraphics[width=0.9\textwidth]{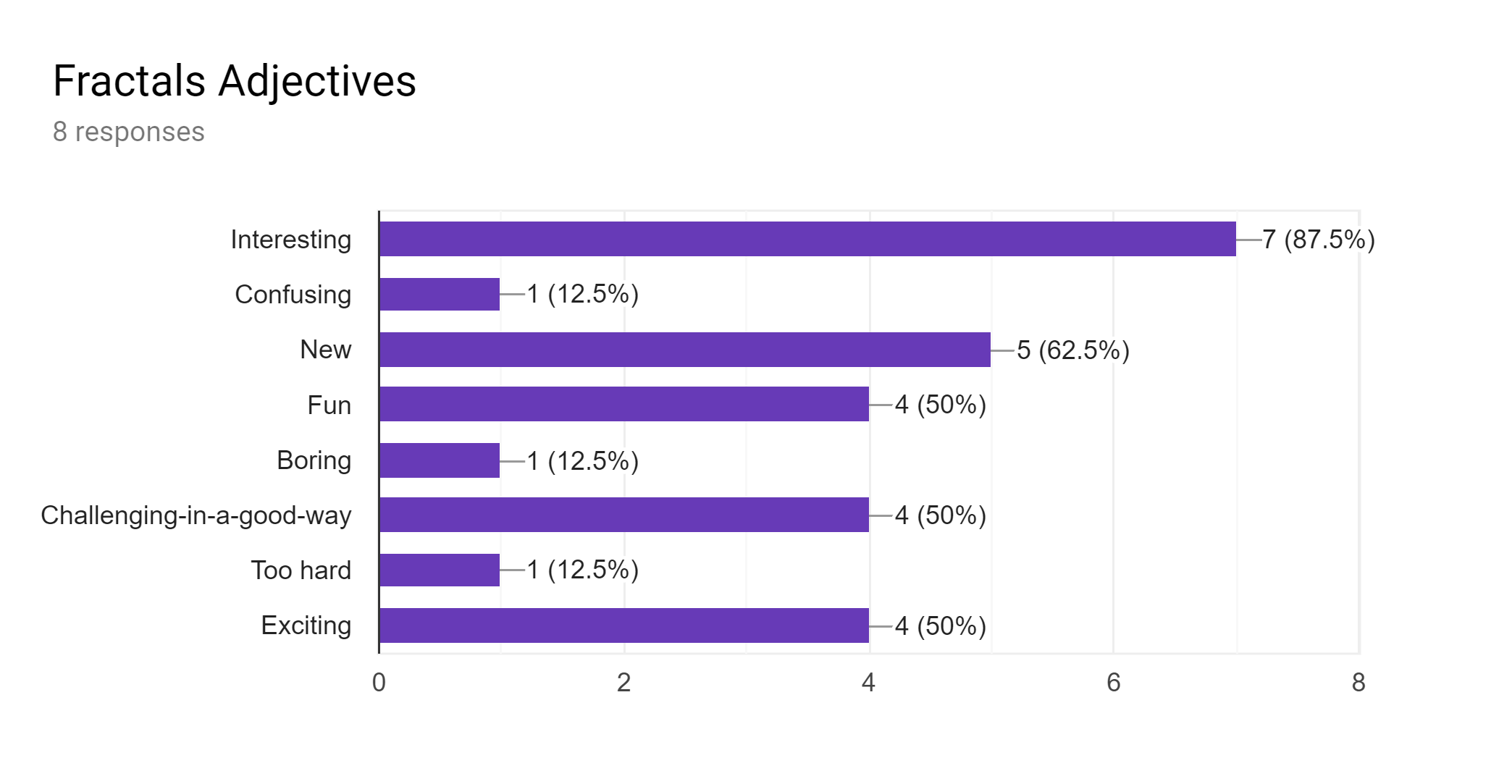}
\end{center}
\caption{Some responses of the students after the second week}
\label{adj2}
\end{figure}

Fig.~\ref{adj2} shows some student feedback after the second week. The main outcome of the course seemed
positive. However, there were some issues that could be improved in a future course.

\begin{enumerate}
\item 
Definition of convergence was introduced quite late. Practical questions arose earlier
than the need of notion of convergence. This had a (probably negative) side effect that students were
unsure what convergence exactly means. Maybe there should be some dedicated time to clarify this
in a future course.

\item The example on the ``monthly interest'' should be eventually done with other numbers to avoid
mixing the values of different parameters in the Bernoulli inequality in general.

\item The triangle inequality was used extensively during studying the Mandelbrot set. Next time
it would be useful to spend some time with the preparation of this very important basic theorem.
Eventually students may have more chance to find their own proofs of Theorems \ref{t1} and \ref{t2}.

\item Also, the proof technique in Theorem \ref{t2} might have been elaborated more. The idea was completely
unknown for all students. With more preparation they could have found their own way to a proof.

\item The ``technique of Day 7'' appeared in several proofs. It was used by some students quite
quickly, but other students found it difficult. Therefore, some more introduction of this method
could be fine.

\item We used some well-known basic theorems that may be not trivial for beginners. Among others,
the statement ``a monotone bounded sequence is convergent'' had a special focus. This was quickly
explained, but maybe more explanation could be even better, with examples (or even by providing a proof).

\item It remained completely hidden why the two definitions of the dragon curve by folding and by ``tiling''
are equivalent. We simply used our intuition. But, in fact, such a precise analysis of the
properties of the dragon curve was not the main purpose of this course. Instead, it was a special call
to observe surprising connections between different approaches of the same problem. Further
investigation of additional questions remained open and expected to think about them
by the students.

\item Tiling the plane by a \textit{high number of pieces of the same iteration} of the dragon curve is actually
different from taking \textit{4 pieces of a high iteration}. Both questions can be answered by
using similar ideas, but a detailed answer should distinguish between them. This kind of difference
was not highlighted during the course.

\item Printing several iterations of the dragon curve remained an open project---the walls of the proposed
3D objects had either a zero width or they were overlapping and therefore unable to print. 
We refer to a related project that proposes an interpolated output---it is available at \cite{dev-frac-cur}.

\end{enumerate}

A final conclusion is that novel technology plays an important role in visualizing higher mathematics---even
for very gifted students. More students were already familiar with programming languages and were
able to write simple programs that show the first few iterations of some fractals. GeoGebra's
computer algebra capabilities and XaoS' real-time visualization helped in finding conjectures very quickly.
Realization of 3D objects that have connections to fractals seemed to be helpful in
understanding some concepts of iterations and self-similarity.

On the other hand, all these instruments were used experimentally. Further research may be needed to
find the most powerful ways to help understanding some basic concepts of fractals for a wider audience
of students as well.

\section{Acknowledgments}
Emily Castner kindly collected and
processed the students' feedback. Keri Celeste, the counselor of the group of children,
gave useful feedback on the students' work. The author's ideas were supported
by several colleagues, including Csaba Biro.

Some useful feedback on the
first versions of this paper was
given by Jameson Parker. Special thanks to Benedek Kov\'acs for ``The Thing''.

\bibliography{kovzol,external}

\end{document}